\documentclass[oneside,a4paper,12pt]{article}
\usepackage{amssymb}
\usepackage{amsmath}
\usepackage{amsthm}
\usepackage[square]{natbib}
\usepackage{texdraw}
\usepackage{epsfig}
\usepackage[active]{srcltx}
\usepackage[all]{xy}
\usepackage{varioref}
\newtheorem{thm}{Theorem}[section]

\newtheorem{cor}[thm]{Corollary}
\newtheorem{lem}[thm]{Lemma}
\newtheorem{prop}[thm]{Proposition}

\newtheorem{defin}[thm]{Definition}

\newtheorem{rem}[thm]{Remark}

\newcommand{\CP}{\mathbb{C\mkern1mu P}}               
\newcommand{\HP}{\mathbb{H\mkern1mu P}}               
\newcommand{\C}{\mathbb{C}}

 \newcommand{\N}{\mathbb{N}}
 \newcommand{\R}{\mathbb{R}}
\newcommand{\Z}{\mathbb{Z}}

\newcommand{\G}{\mathsf{G}}
\renewcommand{\L}{\mathsf{L}}
 
 \newcommand{\T}{\mathsf{T}}

 \newcommand{\SO}{\mathsf{SO}}
 
\newcommand{\Spin}{\mathsf{Spin}} 

\newcommand{\into}{\hookrightarrow}

\renewcommand{\H}{{\mathsf{H}}}
\newcommand{\K}{\mathsf{K}}
\newcommand{\U}{{\mathsf{U}}}
\newcommand{\SU}{{\mathsf{SU}}}
\newcommand{\Sp}{{\mathsf{Sp}}}

\renewcommand{\S}{{\mathsf{S}}}
\newcommand{\diag}{{\text{diag}}}
\renewcommand{\O}{{\mathsf{O}}}

\newcommand{\rk}{{\text{rank}}}
\title{Cohomogeneity one manifolds with positive Euler characteristic}
\author{Philipp Frank, WWU M\"unster}

\begin{document}
\maketitle
\begin{abstract}
  We classify those manifolds of positive Euler characteristic on which
a Lie group $\G$ acts with cohomogeneity one, where $\G$ is classical simple
\end{abstract}

\section*{Introduction}\addcontentsline{toc}{section}{Introduction}

Cohomogeneity one manifolds are Riemannian manifolds with an action of a Lie
group such that the orbit space is one dimensional. They can be seen as
generalisations of homogeneous spaces, but offer a richer structure which has
been used for explicit constructions in the past. In particular, Bergery used
cohomogeneity one manifolds to construct invariant Einstein metrics, and later
Bryant and Salamon found cohomogeneity one metrics with
exceptional holonomy groups, which is not possible in the homogeneous case. For
details and many references see e.g.\ \cite{AA:1992}, \cite{AP:1997},
\cite{AA:1993} and \cite{Alekseevsky:1992}.

An interesting connection arises with nonnegative or positive sectional
curvature. Most constructions of nonnegatively or positively curved manifolds
arise from product and quotient constructions, starting with Lie groups and
their bi-invariant metric. Grove and Ziller found a large class of nonnegatively
curved manifolds within the cohomogeneity one category in
\cite{GroveZiller:2000}: If there are two orbits of codimension two, there is a
metric of nonnegative sectional curvature (this in particular includes the
principal $\L$-bundles over $\S^4$ with $\L=\SO(3)$ or $\L=\SO(4)$). Later Grove
and Ziller showed that every cohomogeneity one manifolds supports a metric of
nonnegative Ricci curvature, and one of positive Ricci curvature if and only if
the fundamental group is finite (\cite{GroveZiller:2002}). The naturally arising
question if every cohomogeneity one manifold with finite fundamental group
allows a metric of positive sectional curvature was answered negatively in
\cite{GVWZ}, where the Brieskorn variety is shown to be a counter example.

This work is concerned with the classification of cohomogeneity one manifolds
with positive Euler characteristic where the acting group $\G$ is classical
simple. The first result is the following theorem, which hints at the importance
of the classification result. The action is called \emph{primitive} if there is
no subgroup $\L\subset\G$ with a $\G$-equivariant map $M\to\G/\L$. Two actions
are called \emph{orbit equivalent} if there is a diffeomorphism of the underlying
manifolds preserving orbits.

\begin{thm}\label{thm:nonsimple}
Suppose the compact, connected Lie group $\G$ acts primitively on the manifold
$M$ with positive Euler characteristic such that $M/\G=[0,1]$. Suppose there is
no normal subgroup of $\G$ that acts orbit equivalently,
and that $\G$ is not simple. Then one of the following applies:

\begin{enumerate}
  \item The action of $\G$ is equivalent to a cohomogeneity one action of a
    rank 1 symmetric space.
  \item $\G$ is covered by $\G'\times \S^3$ where $\G'$ is simple and one of the singular
orbits has codimension 3.
\end{enumerate}
\end{thm}

The major part of this work is the proof of the following
classification:

\begin{thm}\label{thm:main}
Each simply connected, primitive cohomogeneity one $\G$-manifolds of positive
Euler characteristic, where $\G$ is a classical simple group acting almost
effectively,  is, up to equivalence (that is, up to an outer automorphism of
$\G$, up to $\G$-equivariance) one of the following:

\begin{itemize}
    \item[$\bullet$] A linear action on a symmetric space
    \item[$\bullet$] One of the following homogeneous spaces with a linear
	    action: A Grassmanian or one of

    \begin{displaymath}
      \frac{\Sp(n+1)}{\Sp(n-k+1)\U(k)},\quad\frac{\SO(n+1)}{\SO(n-2k+1)\U(k)}
    \end{displaymath}

  \item[$\bullet$] Contained in the the tables given in appendix
	  \ref{subsec:result}.
\end{itemize}
\end{thm}

While the classification tables were moved to its own subsection within the
appendix, it should be noted that the number of examples in each dimension is
finite. We know of no abstract reason for this. We will prove theorem
\ref{thm:nonsimple} in section 2, after we have established the basics and
notation in section 1.  Because of the diagram-chase like qualities of the
classification, we give an overview over the general procedure in section 3. In
sections 4-9, the actual classification is carried out.

This work is part of my PhD thesis which would not have been possible without my
advisor Burkhard Wilking, who I am indebted to for constant support throughout
the years. I'd also like to thank Wolfgang Ziller for his help during my visit
in 2007, and the referees of this paper for valuable advice.

\section{Preliminaries}\label{sec:prelims}

In this section, we will cover the basics of cohomogeneity one manifolds, and
list the properties most usefull to us, especially those connected to the
topology of $M$. We will also treat some properties of representations of
certain simple Lie groups, which we will make use of later in the
classification.

\subsection{Cohomogeneity one actions}\label{subsec:cohom1actions}

Suppose a compact Lie group $\G$ acts on the manifold $M$ such there is
a 1-codimensional orbit $\G/\H$. It is known (see \cite{Br} or \cite{Mostert})
that $M/\G$ is homeomorphic to either $\S^1$ or $[-1,1]$ in the compact case, 
or $\R$ or $[0,1)$ in the noncompact case. The structure of all cases can be
found in the pioneer work of \cite{Mostert}, and also (in more detail)
in various papers by Alekseevsky, see e.g. \cite{AA:1992}. We list the most
important properties:

\subsubsection{The noncompact case.}\label{subsubsec:noncompact} If $M/\G=\R$,
all
orbits are regular and $M$
is actually a fibre bundle over $M/\G$ with fibre $\G/\H$, which is neccessarily
trivial, since $\R$ is contractible. Therefore, $M$ is $\G$-equivalent to
$\G/\H\times \R$. In the case $M/\G=[0,1)$, there is exactly one singular orbit,
say $\G/\K$ (where $\K\supset \H$). Then $\K/\H$ is a sphere $\S^m\subset \R^{m+1}$,
and the action of $\K$ on $\S^m$ can be extended linearly to $\R^{m+1}$. $M$ is
$\G$-equivalent to $\G\times_\K\R^{m+1}$, on which $\G$ acts by left translation on
the first factor.

\subsubsection{The compact case.}\label{subsubsec:compact} If $M/\G=\S^1$, all
orbits are regular and $M$ is a fibre bundle over $\S^1$ with fibre $\G/\H$. It's
easily seen by the homotopy sequence of that bundle that $\pi_1(M)$ is infinite
in this case. Also we have $\chi(M)=0$, so we will not be concerned with this
case.

The most interesting case is $M/\G=[-1,1]$, which is the one we are concerned
with in this paper, and which we therefore lay out in more detail. Choose
any $\G$-invariant Riemannian metric on $M$. There are two singular orbits, so
choose a minimal geodesic $c:[-1,1]\to M$ between them, such that $\pi\circ
c=id_{[-1,1]}$ for the projection $\pi:M\to[-1,1]$. Denote the isotropy groups
by $H=\G_{c(0)}$ and $K_{\pm}=\G_{c(\pm1)}$. By the slice theorem (see e.g.
\cite{Br}) there are tubular neighbourhoods $\G\times_{\K_{\pm}}D_{\pm}$ of the
singular orbits $\K_{\pm}$, where $D_{\pm}$ is the disc of radius 1 normal to
the singular orbit $\K_{\pm}$. Since the action of $\G$ on $M$ has cohomogeneity
one, the action of $\K_{\pm}$ on $\S^{l_\pm}=\partial D_{\pm}$ is transitive,
with isotropy
group $\H$. Now $M$ is constructed from the tubular neighbourhoods by gluing
them along their common boundary $\pi^{-1}(0)\simeq \G/\H$:
\begin{displaymath}
 M\simeq \G\times_{\K_{-}}D_{-} \cup_{\G/\H}\G\times_{\K_{+}}D_{+}
\end{displaymath}
This uniquely determines $M$ in terms of the groups $\G\supset \K_-,\K_+\supset
\H$: The action of $\K_{\pm}$ on the sphere $\S^{l_\pm}=\K_{\pm}/\H$ is linear and
can therefore be extended to an action on the ball $D_{\pm}$, which in turn
allows $M$ be constructed as above. Note that the only condition imposed besides
the inclusions is that $\K_{\pm}/\H$ is a sphere of any dimension. We call the
collection of groups $\{\G,\K_-,\K_+,\H\}$ the \emph{group diagram} of $M$, where it
is understood that $\K_{\pm}\subset \G$ and $\H\subset \K_{\pm}$. The group diagram
is of course not uniquely determined by $M$, since we did not assume
effectiveness of the action of $\G$: If $\G\supset \K_-,\K_+\supset \H$ is the group
diagram of an ineffective cohomogeneity one manifold with ineffective kernel
$\H'$, the effective version of the $\G$ action has group diagram $\G/\H'\supset
\K_-/\H',\K_+/\H\supset \H/\H'$.  But even effective group diagrams are not uniquely
determined by $M$ (e.g., spheres can be represented by many such diagrams).

\subsection{Equivalence, uniqueness of diagrams and primitivity}\label{subsec:uniq}
\emph{For the remainder of this work, we assume that  we are in the second case
described in \ref{subsubsec:compact}, and we use the notation introduced
there.}

If $M_1, M_2$ carry a cohomogeneity one action of $\G$, we will
call $M_1$ and $M_2$ $\G$-\emph{equivariantly diffeomorphic} if there is  a
diffeomorphism $\psi:M_1\to M_2$ such that $\psi(gp)=g\psi(p)$ for all $g\in
\G,p\in M_1$. This will determine the group diagram of a cohomogeneity one
$\G$-manifold $M$ up to the following operations (see \cite{GWZ}):

\begin{itemize} \item Switching $\K_-$ and $\K_+$ \item Conjugating each group of
    the group diagram with an element of $\G$ \item Replacing $\K_-$ by
    $a\K_-a^{-1}$ for some $a\in N(\H)_0$, the unity component of the normalizer
    of $\H$ in $\G$, while keeping $\H$ and $\K_+$ unmodified. \end{itemize}

Of course, the classification works the other way around: We will classify
the possible diagrams up to the listed operations and therefore get a list for
the equivariant diffeomorphism types of cohomogeneity one $\G$-manifolds for a
fixed Lie group $\G$. We will not classify cohomogeneity one manifolds up to
equivariant diffeomorphism, but up to a slightly coarser equivalence
relation:

\begin{defin}
  Two cohomogeneity one $\G$-manifolds $M_1,M_2$ are called \emph{equivalent}, if they are
  equivariantly diffeomorphic up to a automorphism of $\G$, i.e.\ if there is
  an automorphism $\varphi$ of $\G$ and a diffeomorphism $\psi:M_1\to M_2$ such
  that $\psi(gp)=\varphi(g)\psi(p)$ for all $g\in\G,p\in M_1$.
\end{defin}

It is clear that applying an automorphism of $\G$ to all groups in the diagram
of a cohomogeneity one manifold yields an equivalent cohomogeneity one manifold.

From a diagram $\G\supset \K_-,\K_+\supset \H$ and an embedding $\G\into \G'$ we get
a cohomogeneity one $\G'$-manifold from the diagram $\G'\supset \K_-,\K_+\supset
\H$. To further diminish the items of the classification we will require
minimality with respect to this extension process:

\begin{defin}\label{def:prim}{A comohogeneity one $\G$-manifold $M$  is called
\emph{primitive}, if for every group diagram $\G\supset \K_-,\K_+\supset \H$ of
$M$ there is no subgroup $\L\subset \G$ such that $\K_-,\K_+\subset \L$.}
\end{defin}

Nonprimitivity of $M$ is equivalent to the existence of a $\G$-equivariant map
$M\to \G/\L$ for some subgroup $\L\subsetneq \G$. In this case, $M$ is
equivalent to $\G\times_\L N$, where $N$ is the cohomogeneity one manifold
determined by the group diagram $\L\supset \K_-,\K_+\supset \H$ (see
\cite{AA:1992}). $M$ is called an \emph{extension} of $N$. The action of $\L$ on
$N$ does not need to be effective, even if the one of $\G$ on $M$ is.

As a motivation for studying primitive $\G$-manifolds, note that if $N$ has
non-negative curvature, so has $\G\times_\L N$. The technical
value of primitivity lies mostly in the following lemma, taken from \cite{GWZ}:

\begin{lem}\label{lem:prim}
 Assume the primitive cohomogeneity one action of $\G$ on $M$ is effective
(almost effective). Then the intersection $\H_-\cap \H_+$ of the ineffective
kernels $\H_\pm$ of $\K_\pm/\H$ is trivial (finite).
\end{lem}

\subsection{The topology}\label{subsec:fundamentalgroup}
The description of $M$ as a double disc bundle in \ref{subsubsec:compact} of
course allows a computation of the fundamental group in terms of the group
diagram. Vice versa, assuming $M$ being simply connected, one arrives at the
following properties of the groups involved:

\begin{prop}[\cite{GWZ}]\label{prop:spheredim1}
 Let $M$ be a simply connected cohomogeneity one $\G$-manifold. Then there
are no exceptional orbits, and, in the notation of \ref{subsubsec:compact}, we
have $l_\pm\geq 1$, that is, dim $\K_\pm>$ dim $\H$.
\end{prop}

A more direct computation using van Kampen's theorem, carried out in
\cite{Hoelscher:2007}, yields the following usefull properties:

\begin{prop}\label{prop:spheredim2}
 With the assumptions and notation of \ref{prop:spheredim1}:
\begin{itemize}
 \item If $l_+>1$ and $l_-\geq 1$, then $\K_-$ is connected
 \item If both $l_\pm >1$, then all of $\H,\K_-,\K_+$ are connected
\end{itemize}
\end{prop}

We will utilize this proposition frequently without further mention: The only
way the group $\K_+$ can be non-connected is when the ``other sphere'' $\K_-/\H$
is 1-dimensional.

Another result concerns the Euler characteristic of $M$ (\cite{AP:1997}). Since
$\G\times_{\K_\pm}D_\pm$ contains the singular orbit $\G/\K_\pm$ as a deformation
rectract, the Mayer-Vietoris exact sequence can be applied to obtain the
following exact sequence:
\begin{displaymath}
 \ldots \to H^i(M) \to H^i(\K_-)\oplus H^i(\K_+)\to H^i(\G/\H)\to
H^{i+1}(M)\to\ldots
\end{displaymath}

This allows the following corollary:

\begin{cor}\label{cor:Eulerchar}
 Let $M$ be a cohomogeneity one $\G$-manifold with orbit space $[-1,1]$. Using
notation from \ref{subsubsec:compact}, the Euler characteristic of $M$ is given
by
\begin{displaymath}
 \chi(M)=\chi(\G/\K_-)+\chi(\G/\K_+)-\chi(\G/\H)
\end{displaymath}
In particular,  $\chi(M)>0$ implies that one of $\chi(\G/\K_\pm)$ is greater
than zero, which is equivalent to one of the singular isotropy groups
$\K_\pm$ having maximal rank in $\G$.
\end{cor}

Given the cohomogeneity one diagram, it is easy to compute the Euler
characteristic (see e.g.\ \cite{Wang}). For a homogeneous space $\G/\K$ where
$\G,\K$ are connected, compact with the same rank, we have

\begin{displaymath}
	\chi(\G/\K)=o(\G)/o(\K) 
\end{displaymath}

where $o$ denotes the order of the weyl group. For $\K$ not connected we only
need to note that $\G/\K$ is a finite cover of $\G/(\K)_0$.

For $\G$ simple, the maximal subgroups of maximal rank have been
classified by Borel and Siebenthal (see e.g. \cite{GoGr}), up to an
automorphism of $\G$. This can easily be
generalized to semisimple and compact groups, as well as to the non-maximal
subgroups. We list the result for the classical groups in table
\ref{table:maxranksubgroups}.

\begin{cor}\label{cor:corH}
  Suppose $\chi(M)>0$ and that $\K_+$ has the same rank as $\G$. Then the group
  $\H$ has corank 1 in $\G$, i.e.\ the sphere $\S^{l_+}$ is odd dimensional.
\end{cor}

\begin{proof}
  Suppose $p\in \G/\K_+$ is in the singular orbit given by $\K_+$. Then the
  dimension of the slice representation at $p$ is given by $\dim M - \dim
  \G/\K_+$. But both $M$ and $\G/\K_+$ have positive Euler characteristic, so
  they are even dimensional, and so is the slice representation at $p$.
  Therefore the sphere of the slice at $p$ is odd dimensional.
\end{proof}

We note that in the case that $\text{corank}(\H)=1$ and $\text{corank}(\K_+)=0$
we can derive $\chi(M)>0$. This is not true in general for
$\text{corank}(\H)=0$.
\subsection{Properties of the isotropy groups} Again, we carry over the
assumptions and notation from \ref{subsubsec:compact} and \ref{lem:prim}.
Since $\G$ is compact, it is finitely covered by a group $\tilde \G$ of the form
$\T^k\times \G_1\times\cdots\times \G_l$ where $\T^k$ is a $k$-dimensional torus,
the center of $\tilde \G$, and $\G_1,\ldots, \G_l$ are simple normal subgroups. We
call $k+l$ the \emph{number of factors of $\G$}, which is well-defined since $\G$ is
a finite quotient of $\tilde \G$. The fact that if $\G$ acts transitively on a
sphere, the isotropy group has at least $k+l-1$ factors (see table
\ref{table:sphereactions}) leads to the following lemma:

\begin{lem}\label{lem:numberoffactors}
 If $\G$ acts almost effectively and primitively, $\K_+$ has at most 4 factors.
\end{lem}
\begin{proof}
 Denote the number of factors of a group $\H$ by $f(\H)$. Since $\H_\pm$ are
 normal subgroups of $\H$, it is clear that

 \begin{displaymath}
   f(\H_+)+f(\H_-)-f(\H_+\cap\H_-)\leq f(\H)
 \end{displaymath}

 so that

 \begin{displaymath}
   f(\H_+)+f(\H_{-})-f(\H)\leq f(\H_+\cap\H_-)
 \end{displaymath}

 By the classification of effictive transitive actions on spheres (see table
 \ref{table:sphereactions}) we have $f(\H)-f(\H_\pm)=f(\H/\H_\pm)\leq 2$, which implies 

 \begin{displaymath}
   2f(\H)-f(\H_-)-f(\H_+)\leq 4
 \end{displaymath}

 and consequently

 \begin{displaymath}
   f(\H)-4\leq f(\H_+)+f(\H_-)-f(\H)\leq f(\H_+\cap\H_-)
 \end{displaymath}

 But now $f(\H_+\cap\H_-)=0$ by lemma \ref{lem:prim}, so that $f(\H)\leq 4$. By
 the classification of sphere actions we know $f(\K_+)\leq 5$, but if
 $f(\K_+)=5$, we know that $\K_+/\H=\S^1$ or $\K_+/\H=\S^3$, and
 $(\H)_0=(\H_+)_0$ has 4 factors. This is a contradiction to
 $\H_+\subset\H/\H_-$, given by primitivity. So $f(\K_+)\leq 4$.
 
\end{proof}

\begin{rem}\label{rem:numberoffactors}
 Since we are concerned with almost effective actions, any factor of $\K_+$ that
acts trivially on $\K_+/\H$ (that is, is contained in $\H_+$), cannot act
trivially on $\K_-/\H$, so it has to occur as a factor of an isotropy group of an
effective transitive sphere action in table \ref{table:sphereactions}, up to
finite quotients. If $\K_+$ has 4 factors $\K_1\cdots \K_4$, we can after
rearranging the order assume that $\K_1\K_2$ act transitively almost effectively
on a sphere, and $\K_3\K_4$ appear as the isotropy group of a transitive
effective action on  a sphere, up to finite quotients. This
will limit the choices when looking  for possibilities for $\K_+$ from table
\ref{table:maxranksubgroups}.
\end{rem}

\begin{lem} The kernel of the action is the largest normal subgroup shared by
  $\G$ and $\H$ 
\end{lem} 

\begin{rem}
  \label{rem:finitenormalsubgroups}
  For $\G$ simple, the only normal subgroups are finite and therefore central. So
  the kernel of the action is the intersection of $\H$ with the center of $\G$.
\end{rem}

\subsection{Known classification results}\label{subsec:known}

\subsubsection{Cohomogeneity one manifolds of low
dimension}\label{subsubsec:lowdim}

Cohomogeneity one manifolds of dimension up to seven have been classified, with
no assumption on the Euler characteristic: Neuman (\cite{Neumann}) classified
those of dimension three, Parker (\cite{Parker}) those of dimension four (with
one omission, as observed by C. Hoelscher), and Hoelscher
(\cite{Hoelscher:2007}) classified those of dimension five to seven. The only
examples of positive euler characteristic are symmetric spaces, even if the
group $\G$ is not assumed to be simple.

\subsubsection{Cohomogeneity one manifolds with a fixed
point}\label{subsubsec:fixedpoint}

If the action of $\G$ on $M$ has a fixed point, the classification is
particularly easy: There are the obvious actions with two fixed points on
spheres, and the following groups acting on compact rank one symmetric spaces:
 \begin{align*}
 \CP^n&: \SU(n),\U(n)\\
 \HP^n&: \Sp(n),\Sp(n)\times\Sp(1),\Sp(n)\times\U(1)  \\
\CP^{2n+1}&:\Sp(n), \Sp(n)\times \U(1)\\
Ca\mathbb P^2&:\Spin(9)
 \end{align*}
The details can be found in \cite{Hoelscher:2007}

\subsubsection{Cohomogeneity one manifolds with positive euler characteristic}
The pioneer work on this subject was done in the already cited paper by
\cite{AP:1997}. The present work enhances its result in the following ways.
Firstly Alekseevsky and Podest\'a only give a classification under the
assumption of $\G$ either having a fixed point or at least one of the spheres
$\K_\pm/\H$ being effective. Secondly it is desireable not only give the
isomorphism classes of the lie algebras of the lie groups involved, but to give
embeddings and components as well. As an example, there are 4
non-equivalent cohomogeneity one $\SO(2n)$ manifolds where $\frak
{k}_\pm\simeq\mathbb R\oplus\frak{su}(n)$ and $\frak h\simeq\mathbb
R\oplus\frak{su}(n-1)$. Note that all of $\G,\K_\pm,\H$ are connected in all of
those 4 examples, while 2 of them are primitive. As a second example, there are
two $\SU(3)$-cohomogeneity one manifolds where $\frak k_+\simeq
\frak{so}(3),\frak k_-\simeq \mathbb R^2,\frak h\simeq\R$, both of which are
primitive, and for one of those 2 examples, $\K_-$ and $\H$ are not connected.

Alekseevsky and Podest\'a do not require $\G$ to be simple.
\section{Proof of theorem \ref{thm:nonsimple}}
We start with a little lemma concerning orbit equivalent subactions, which will
be used in the proof:

\begin{lem}
Let $M$ be the cohomogeneity one $\G$-manifold given by the diagram
$\H\subset\K_-,\K_+$, and suppose $\G=\G_1\times\G_2$. If the projection of $\H$
onto the second factor is $\G_2$, then the subaction of $\G_1\times 1$ on $M$ is
also cohomogeneity one.
\end{lem}
\begin{proof}
The claim can be tested on any orbit $\G/\G_x$, where $\G_x$ is one of
$\H,\K_\pm$.  For any $(g_1,g_2)\G_x\in\G/\G_x$, there is some $(h_1,g_2)\in\H$
by the assumption. Then $(g_1,g_2)\G_x=( g_1 h_1^{-1},1)( h_1, g_2)\G_x=( {g_1
h_1^{-1}},1)\G_x$, so $( g_1, g_2)\G_x$ is in the $\G_1\times 1$-orbit of $( 1,
1)\G_x$.
\end{proof}

We can now prove theorem \ref{thm:nonsimple}.

\begin{proof}
By virtue of $M$ having positive Euler characteristic, we can assume $\K_+$ has
maximal rank in $\G$. We also assume $\G=\G_1\times\dots\times\G_l$ where each
$\G_i$ is either simple or $\S^1$. Then also
$\K_+=\K_+^1\times\dots\times\K_+^l$, where each $\K_+^i$ is a subgroup of
$\G_i$ of maximal rank, and we can assume that $\K_+^1$ acts transitively on the
sphere $\K_+/\H$. Let $pr_i$ be the projection from $\G$ onto the $i$-th factor.
By the lemma above, $pr_i(\H)\ne\G_i$ for all $i=1,\dots,l$.

We claim $l=2$. Suppose $l>2$ and let $p_2:\G\to\G_2\times\dots\times\G_l$ be the projection.
Because $\K_+^1$ acts transitively, we have $p_2(\K_+)=p_2(\H)$, and primitivity
implies $p_2(\K_-)=\G_2\times\dots\times\G_l$ (otherwise, $\K_\pm\subset
p_2^{-1}(p_2(\K_-))\ne\G$, a contradiction to primitivity). But
$pr_2(\H)\ne\G_2$, so $p_2^{-1}(\G_2)\cap\K_-$ acts transitively on $\K_-/\H$,
which implies $p_3(\K_-) = p_3(\H)$ for the projection
$p_3:\G\to\G_3\times\dots\times\G_l$. By the lemma above, $\G_1\times\G_2$ acts
orbit equivalent, a contradiction. Therefore $l=2$.

For now suppose that $\K_-\cap\G_2$ is not finite. We have $pr_2(\K_-)=\G_2$, so
$\K_-\cap\G_2$ is a normal subgroup of $\G_2$ (if $(1,k_2)\in\K_-\cap\G_2$ and
$g_2\in\G_2$ is arbitrary, then there is some $(k,g_2)\in\K_-$, and
$(1,g_2)(1,k_2)(1,g_2^{-1})=(k,g_2)(1,k_2)(k^{-1},g_2^{-1})\in\K_-\cap\G_2$). By
assumption, $\K_-\cap\G_2=\G_2$. But $pr_2(\H)\ne\G_2$, so $\G_2\subset\K_-$
acts transitively on $\K_-/\H$, and therefore $pr_1(\K_-)=pr_1(\H)$. By
primitivity, $pr_1(\K_+)=\G_1$, and from $\K_+$ having the same rank as $\G$ we
can deduce that $\G_1$ is a normal subgroup of $\K_+$. We now divide cases by
the dimension of the sphere $\S^{l_-}$:

\begin{itemize}
  \item[$\bullet$] Suppose $l_-$ is even, that is, the rank of $\K_-$ is the
    same as the rank of $\H$. That implies that $pr_2(\H)$ has full rank in
    $\G_2$ and therefore $\H$ is a product subgroup of $\G$. So we have 

    \begin{align*}
	\K_-&=\K_-^1\times\G_2\\
	\H&=\K_-^1\times\H_2\\
	\K_+&=\G_1\times\H_2
    \end{align*}
    where $\G_1/\K_-^1$ and $\G_2/\H_2$ are spheres. This is easily recognized
    as a so called sum-action on a sphere of even dimension (see e.g.\ \cite{Hoelscher:2007}).
  \item[$\bullet$] If $l_-$ is odd, that is, the rank of $\K_-$ is
    $\rk(\H)+1=\rk(\G)$, then we have the following situation: $
    \K_-=\K_-^1\times\G_2, \K_+=\G_1\times\K_+^2$. If we define
    $\H_i:=\H\cap\G_i$, we have that $\G_i/\H_i$ is a sphere for $i=1,2$, and
    $\H_1\times\H_2$ is a normal subgroup of $\H$ of corank 1. Then
    we can find a rank 1 normal subgroup $\Delta\H$ of $\H$ that commutes with
    $\H_1\times\H_2$ such that $\H=(\H_1\times\H_2)\Delta\H$. We have
    $\K_-^1=pr_1(\H)$ and $\K_+^2=pr_2(\H)$.

    Now consider the cohomogeneity one $\G_1\times\G_2$-manifold given by the
    group diagram $\H_1\times\H_2\subset\H_1\times\G_2,\G_1\times\H_2$, which is
    a sphere $\S^{2n+1}$ of odd dimension as above. We claim that $M$ is the
    quotien of $\S^{2n+1}$ by a free action of $\H/(\H_1\times\H_2)$. For that we
    only need to consider the following actions of $\H/(\H_1\times\H_2)$ on the
    orbits of $\S^{2n+1}$:

    \begin{itemize}
      \item $\H/(\H_1\times\H_2)$ acts freely on $\G_1/\H_2\times\G_2/\H_1$ with
	quotient $\G_1\times\G_2/\H$.
      \item $\H/(\H_1\times\H_2)$ acts freely on
	$\G_1/\G_1\times\G_2/\H_2=\G_2/\H_2$ with quotient $\G_2/pr_2(\H)$.
      \item $\H/(\H_1\times\H_2)$ acts freely on
	$\G_1/\H_1\times\G_2/\G_2=\G_1/\H_1$ with quotient $\G_1/pr_1(\H)$.
    \end{itemize}

    In conclusion, the orbits of the action of $\G$ on the quotient of
    $\S^{2n+1}$ by $\H/(\H_1\times\H_2)$ are those of the action of $\G$ on $M$,
    which finishes this part of the proof (see section \ref{sec:prelims}).

    Lastly, consider $\K_-\cap\G_2$ finite. Then
    $\rk(\G_2)=1$ and $\K_-/\H$ is even dimensional (if
    $\rk(\K_-)=\rk(\G)$, i.e.\ the sphere is odd dimensional, then $\K_-$ is a
    product with $pr_2(\K_-)=\G_2$, a contradiction; also $\K_-\cap
    \G_2$ has corank 0 or 1 in $\G_2$, which implies $\rk(\G_2)\leq 1$). We also
    have $pr_2(\K_-)=\G_2$ and $pr_2(\H)\ne \G_2$, so $\G_2$ is covered by
    $\S^3$ and $\K_-/\H=\S^2$.
\end{itemize}
\end{proof}

\section{The general procedure}\label{sec:procedure}

In this section we will describe the actual procedure for the classification.
We will treat each of the simple groups separately, and make extra sections for
$\SU(3)$ and $\SU(4)$.

For each classical simple group $\G$, we will first list the result, the table
of group diagrams of cohomogeneity one $\G$-manifolds up to equivalence. In
order to prove this result, we will then list the possibilities for the subgroup
$\K_+$ of maximal rank, combining table \ref{table:maxranksubgroups}, lemma
\ref{lem:numberoffactors} and remark \ref{rem:numberoffactors}. This is simply
an exercise in book-keeping. The conditions listed ensure that the different
cases really are disjoint. Note that for the sake of organisation we divided
cases such as $\SO(n),n\geq 2$ into the cases $\SO(2)$ and $\SO(n),n\geq 3$.

By what was said in section \ref{subsec:uniq}, we can assume $\K_+$ has the
standard block structure. We will then use table \ref{table:sphereactions} to
list the possibilities for the isotropy group $\H$ of the action of $\K_+$ on
the sphere $\K_+/\H$. Again, we can conjugate the diagram by an element of
$\K_+$ to ensure $\H$ is of a given form. After that, we can again use the same
table to list the possibilities for $\K_-$. Lemma \ref{lem:prim} will be used
without further mention to discard some of the possibilities. The last step is
to check the possible embeddings of $\K_-$ into $\G$, i.e.\ which are equivalent
and which give a primitive diagram.

As far as possible, we will give the homogeneous representation for those
examples that are actually homogeneous spaces. Most of the claims are easily
checked, the more involved examples can be found in \cite{GWZ}.

\subsection{The $\Spin$-groups}\label{subsec:procspin}

The cases of the $\Spin$-groups will be divided into two different cases each.
First, we will classify the non-effective actions, i.e.\ those that are actually
action of the  special orthogonal group. After that, we will classify the effective
actions of the $\Spin$-groups. The procedure for the latter ones differs slightly
from the general procedure described above: We will apply the projection
$\pi:\Spin(n)\to\SO(n)$ to the whole diagram and classify the resulting
diagrams. The list of subgroups of maximal rank is easily deduced from table
\ref{table:maxranksubgroups}.

Since we know the action is effective, we have $-1\notin\H$ (where $-1$ is the
element that projects to the identity of $\SO(n)$ but which is not the identity
element of $\Spin(n)$) by remark \ref{rem:finitenormalsubgroups}. This shortens
the list of possibilities for $\K_+$, because it implies that $\H$ does not
contain a subgroup of type $\Spin(n),n\geq 3$.

We have $-1\in\K_+$ from the following fact: The preimage $\pi^{-1}(\K)$ of a
subgroup $\K\subset\SO(n)$ is connected if and only if the inclusion
$\K\into\SO(n)$ induces a surjection on the fundamental group.  This is the case
for all maximal rank subgroups of $\SO(n)$, so their preimages contain $-1$ in
the unity component. This implies that $\pi(\K_+)/\pi(\H)$ is a real projective
space, and the possibilities for $\pi(\H)$ can be deduced from table
\ref{table:sphereactions}. We cannot assume $-1\in\K_-$, so we list
the possibilities for $\pi(\K_-)$ under the assumption that $\pi(\K_-)/\pi(\H)$
is either a sphere or a projective space. For a subgroup of $\SO(n)$ of a given
isomorphism type it is easy to decide whether the above criterion applies, so we
can decide if $\K_-/\H$ is a sphere, and carry on as above.

For convenience of notation, we will always give the $\SO(n)$-diagram
$\pi(\H)\subset\pi(\K_-),\pi(\K_+)$. By abuse of notation, we will discard the
$\pi$, which will not lead to confusion, since everything is discussed in
$\SO(n)$ anyways.

There is another simple fact we will make use of. In the situation described
above, note that $\SO(n)/\pi(\H)$ is not simply connected, because $\Spin(n)/\H$
is a nontrivial cover. This implies in particular that $\pi(\H)\into\SO(n)$ does
not induce a surjection on the fundamental group.

%
%
%

\section{$\G=\SU(3)$}\label{sec:su3}

We claim that up to equivalence, the diagrams of the simply connected primitive
cohomogeneity one $\SU(3)$-manifolds with positive Euler characteristic are
given by table \ref{table:su3}.

\renewcommand\arraystretch{1.4}

\begin{table}[h]\caption{$\SU(3)$-cohomogeneity one manifolds}\label{table:su3}
  \begin{center}
  \begin{tabular}{|l|}\hline 
    $\S^1\subset\SU(2),\U(2)$\\\hline
    $\S^1\subset\S(\U(2)\U(1)),\S(\U(1)\U(2))$\\\hline
    $\S^1\subset\SO(3),\S(\U(1)\U(2))$\\\hline
    $\S^1\subset\SO(3),\T^2$\\\hline
    $\Z_3\SO(2)\subset\Z_3\SO(3),\T^2$\\\hline
  \end{tabular} 
  \end{center}
 \end{table} 
 
\renewcommand\arraystretch{1}

    By the classification of Borel an Siebenthal (see table
    \ref{table:maxranksubgroups}), we have $\K_+=\U(2)$ or $\K_+=\T^2$. If
    $\dim \H\geq 2$, we have $\dim M=\dim \G/\H +1 \leq 7$, so $M$ appears in
    \cite{Hoelscher:2007}. So we will assume $(\H)_0=\S^1$ for the rest of this
    section.
    
    First assume $\K_+=\U(2)$, where $\SU(2)$ is embedded in the lower right
    block. By what was said above and what follows from \ref{subsubsec:un1spheres}, we
    have $\H=\S^1_k$ where

    \begin{displaymath}
      \S^1_k=\{\diag(\bar z,z^{k+1},z^{-k})\mid z\in\S^1\}
    \end{displaymath}

    If $\K_-$ acts almost effectively $\S^{l_-}$, we have
    $\K_-\in\{\U(2),\SO(3),\SU(2)\}$. Otherwise $\K_-=\T^2$, which we will treat
    later in this section, or $\K_-=\S^1_k\SU(2)$ where $\S^1_k$ is normal in
    $\K_-$ and does not intersect $\SU(2)$. The first is only possible for
    $k=-2$ or $k=1$, but in both cases the (unique) $\SU(2)$ in its normalizer
    is intersected.

    For an almost effective action of $\SU(2)$, we need to find an $\SU(2)$ that
    contains $\S^1_k$, which implies $k=-1$ or $k=0$. Both cases are equivalent
    by a change of the last 2 coordinates, which fixes $\K_+$, so we can assume
    $k=0$, which gives  a primitive example, because $N(\H)_0=\T^2$. This is
    $\SU(3)$ acting on $\mathbb HP^2$. If $\K_-\simeq \U(2)$, we will argue that
    we can assume both $\K_\pm$ contain the same maximal torus. Of course, we
    can conjugate the maximal torus of $\K_-$ into the standard one, and by
    changing the conjugation with an element of $\K_+$, we may assume it
    preserves $\S^1_k$. But then both $\S^1_k$ and its conjugate are diagonal,
    so by changing the conjugation with an element of the Weyl group, we may
    assume it actually preserves $\S^1_k$ pointwise, so it is contained in its
    centralizer, which is in $N(\S^1_k)_0$.  Now both $\K_\pm$ contain the same
    maximal torus, so they are conjugate by an element of the Weyl group, and
    since one of those elements fixes $\K_+$, we can assume
    $\K_-=\S(\U(2)\U(1))$. Checking the possible isotropy groups of $\K_-$, we
    see $k=\pm1$. For $k=1$, the normalizer of $\S^1_1$ contains an $\SU(2)$ in
    which we can realize the exchange of the first and the third coordinate,
    transforming $\K_-$ into $\K_+$, so this is not a primitive example. For
    $k=-1$, we have $N(S^1_{-1})_0=\T^2$, so this example is primitive, the 
    Grassmanian $\SU(4)/\S(\U(2)\U(2))$.
    
    Now we will show that there is one example for $\K_-=\SO(3)$. For that we
    will argue that $\K_-$ is conjugate to the standard $\SO(3)$ in
    $N(\S^1_k)_0$. Since there is only one 3-dimensional representation of
    $\SO(3)$, we know $\K_-$ is conjugate to the standard subgroup. But all
    $\SO(2)\subset\SO(3)$ are conjugate, so we may assume this conjugation
    preserves the standard $\SO(2)$, i.e.\ it is in its normalizer
    $\O(2)\times\S^1$, where $\S^1$ is given by $\diag(z,z,\bar z^2)$. But every
    conjugation of $\O(2)$ on $\SO(2)$ can be realized in $\SO(2)$ itself, so we
    can assume the conjugation comes from $N(\SO(2))_0$. But now $\S^1_k$ can
    only be contained in $\SO(3)$ if it's conjugate to $\SO(2)$, so the
    embedding $\S^1_k\into\SU(3)$ must have a 1-dimensional trivial
    subrepresentation, implying $k=0,-1$. Both of those are conjugate to the
    standard $\SO(2)$, so by the argument above we can assume $\K_-$ is the
    standard $\SO(3)$. Since $\S^1_{-1}$ and $\S^1_0$ can be transformed into
    each other by an outer automorphism of $\K_+$ (complex conjugation, which is
    also an automorphism of $\SU(3)$) that leaves $\SO(3)$ invariant, we just
    get 1 example from this case.

    The second and last case to consider is $\K_+=\T^2$. We have $\H_0=\S^1$,
    therefore $(\K_-)_0$ is one of $\SO(3),\U(2),\SU(2),\T^2$.

    If $(\K_-)_0=\T^2$, it is contained in the centralizer of $\H_0$ as well as
    $(\K_+)_0$, so both are conjugate in $N(\H)_0$ in particular. Since a maximal
    torus has  finite index in its normalizer, this shows that no primitive
    example arises in this case.

    If $(\K_-)_0=\SU(2)$, then $\H_0$ is a maximal torus in $\K_-$. By conjugating
    the diagram we may assume $(\K_-)_0$ is given by the lower 2x2-block, and
    $\H_0$ the standard maximal torus therein. This determines the maximal torus
    $\T^2$ in $\SU(3)$, and both of $\K_-$ and $\K_+$ are contained in $\U(2)$.

    For $(\K_-)_0=\U(2)$: If $\H_0$ is regular, i.e. $Z(\H)_0=\T^2\subset\U(2)$,
    this is obviously not primitive. If $\H_0$ is not regular, its isotropy
    representation has 2 equal eigenvalues, and we may assume we have  

    \begin{displaymath}
      (\H)_0=\{\diag(z,z,\bar z^2)\mid z\in\S^1\}
    \end{displaymath}

    Its isotropy representation therefore has 2
    equivalent 2-dimensional factors, which are in fact equivalent in $N(\H)_0$,
    so we can conjugate $(\K_-)_0$ into

    \begin{displaymath} \U(2)=\left\{\begin{pmatrix} \det \bar A & 0\\ 0 & A\\
    \end{pmatrix}\mid A\in\U(2)\right\} \end{displaymath} without changing
    $\T^2$, so there is also no new primitive example.

    In the case $(\K_-)_0=\SO(3)$, $\K_-$ is actually given by the standard
    embedding, for $\SO(3)$ has no outer automorphisms and its only faithfull
    3-dimensional representation is irreducible. The latter also implies that
    its centralizer is given by $\Z_3$, the set of diagonal matrices, so we have
    $N(\SO(3))=\SO(3)\Z_3$.  Since $T^2$ is uniquely determined by $\H_0=\S^1$,
    we obtain two new examples, both of which are primitive for the isotropy
    representation of $\SO(3)$ in $\SU(3)$ is irreducible, so $\SO(3)$ is a
    maximal subgroup not containing $\K_+$. The examples are
    $\S^1\subset\{\SO(3),\T^2\}$ and $\S^1\Z_3\subset\{\SO(3)\Z_3,\T^2\}$. Note
    that $\K_+$ is connected, since $l_-=2$.

\section{$\G=\SU(4)$}\label{sec:su4}

We claim that up to equivalence the diagrams of the simply connected primitive
cohomogeneity one $\SU(4)$-manifolds with positive Euler characteristic are
given by tables \ref{table:so6} and \ref{table:su4}.

\renewcommand\arraystretch{1.4}

\begin{table}[h]\caption{$\SO(6)$-cohomogeneity one manifolds}\label{table:so6}
 \begin{center} \begin{tabular}{|c|}\hline
	$\SO(4)\subset\SO(5),\SO(2)\SO(4)$\\\hline
	$\Z_2\SO(4)\subset\Z_2\SO(5),\SO(2)\SO(4)$\\\hline
	$\SO(2)\SO(3)\subset\SO(3)\SO(3),\SO(2)\SO(4)$\\\hline
	$\SO(2)\SO(2)\subset\SO(2)\SO(3),\U(2)\SO(2)$\\\hline
	$\U(2)\subset\SO(4),\U(3)$\\\hline
	$\T^2\subset\SO(3)\SO(2),\SO(2)\U(2)$\\
	where $\T^2=\{\diag(z_1,1,1,z_2)\}$\\\hline
	\end{tabular}\end{center}
 \end{table} 
 
 \begin{table}[h]\caption{$\SU(4)$-cohomogeneity one manifolds}\label{table:su4}
  \begin{center} \begin{tabular}{|c|}\hline
    $\S^1\SU(2)\subset\S(\U(2)\U(2)),\S(\U(1)\U(3))$\\
    where $\S^1=\{\diag(\bar z^2,z^4,\bar z,\bar z)\}\subset N(\SU(2))$\\\hline
    $\S^1\SU(2)\subset\S(\U(2)\U(2)),\S(\U(1)\U(3))$\\
    where $\S^1=\{\diag(\bar z^2,1,z,z)\}\subset N(\SU(2))$\\\hline
    $\S^1\subset\sigma(\S(\U(1)\U(3))),\S(\U(1)\U(3))$\\
    where  $\S^1=\{\diag(\bar z,z,1,1)\}$ and $\sigma$ exchanges the first two
    coordinates\\\hline
    \end{tabular}\end{center}      
 \end{table} 
 
\renewcommand\arraystretch{1}
    By the classification of Borel and Siebenthal, given in table
    \ref{table:maxranksubgroups}, and remark \ref{rem:numberoffactors}, we
    know $\K_+$ is one of $\U(3), \S^1\SU(2)\SU(2), \S^1\U(2)$. 

    First, we deal with that case that $\Sp(2)$ is contained in any of the
    regular isotropy groups. Since $\rk(\Sp(2))=2<\rk(\SU(4))$, and $\Sp(2)$ is
    a maximal connected subgroup of $\SU(4)$, we can deduce $(\K_-)_0=\Sp(2)$.
    Because $\K_+$ has maximal rank in $\SU(4)$, and the rank of $\K_+$ and $\H$
    can differ by at most 1, we see $\rk (\H)=\rk (\K_-)$, so $\K_-/\H$ must be
    an even dimensional sphere. By the classification of transitive effective
    actions on spheres (see table \ref{table:sphereactions}), we know
    $\H=\Sp(1)\Sp(1)$, where the common central element of $\H,\K_-,\SU(4)$ is
    in the kernel of the action of $\G$, so this is actually an action of
    $\SO(6)=\SU(4)/\{\pm Id\}$, and $\left( \K_- \right)_0=\SO(5),\left(
    \H \right)_0=\SO(4)$. Since $\rk(\K_+)=3$, we have $\K_+=\S^1\SO(4)$. Both
    $\K_+$ and $\H$ can have at most 2 components, so $\S^1$ can act with weight
    1 or 2 on the slice, giving two primitive examples (weight 1 leading to the 
    Grassmanian $\SO(7)/\SO(2)\SO(5)$).

    Now we divide cases by $\K_+$, under the assumption that $\K_-$ does not
    contain $\Sp(2)$ as a factor (which is true for $\K_+$ by the classification
    anyways).

    So now assume $\K_+=\U(3)$, where $\SU(3)$ is the lower right block. If
    $\SU(3)\subset \H$, then it's easily seen that $\K_-=\SU(4)$ by table
    \ref{table:sphereactions}, which is listed in subsection
    \ref{subsubsec:fixedpoint}. So we can assume $\SU(3)\not\subset\H$, which
    implies that $\U(3)$ acts almost effectively on the slice, giving
    $\H=\S^1_k\SU(2)$, where

    \begin{displaymath}
      \S^1_k=\{\diag(\bar z^2,z^{2(k+1)},\bar z^k,\bar z^k)\mid z\in\S^1\}
    \end{displaymath}

    and $\SU(2)$ is the lower right block (see section \ref{subsubsec:un1spheres}).
    We now divide cases by the rank of $\K_-$ and its dimension:

    \begin{itemize}
      \item If $\rk(\K_-)=\rk(\H)$, we have $\K_-/\H=\mathbb S^2$
	and $\K_-/\H_-=\SO(3)$. Moreover, the semisimple part of $\H$ is
	contained in $\H_-$, and therefore $\K_-\subset
	N(\SU(2))=\S^1\SU(2)\SU(2)$. That implies $\K_-=\SU(2)\SU(2)$, and since
	$\H\subset\SU(2)\SU(2)$, we have $k=0$. This leaves one primitive
	example, $\SU(4)$ acting on $\mathbb HP^3$. 
    \end{itemize}

    From now on, we can assume $\rk(\K_-) > \rk (\H)$. We divide cases by the
    dimension of $\K_-/\H$, which is easily seen to be bounded by 5.

    \begin{itemize}
      \item If $\K_-/\H=\mathbb S^1$, we have $\T^2\SU(2)=\K_-\subset N(\H)_0$.
	If $k\ne 2$, we have $N(\H)_0=\T^2\SU(2)\subset\K_+$, so there will be
	no primitive example. For $k=-2$, we have $N(\H)_0=\S^1\SU(2)\SU(2)$,
	and it is easily seen that up to conjugation in $N(\H)_0$ we have
	$\K_-\subset\K_+$. As before, this contradicts primitivity.
      \item If $\K_-/\H=\mathbb S^3$, then $\K_-=\S(\U(2)\U(2))$. Since
	$\K_-/\H$ is a sphere, we have that $\SU(2)\cap\H$ is trivial, where
	$\SU(2)$ is the upper left block. This leads to $k=\pm 1$ and gives 2
	examples, which are obviously primitive and a Grassmanian.
      \item Lastly, consider $\K_-/\H=\mathbb S^5$. This
	implies $\K_-=\U(3)$, and by studying the isotropy representation of
	$\SU(4)/\H$, there are two possibilities for $\K_{-}$, but since
	primitivity implies $\K_-\ne\K_+$, we know $\K_-=\sigma\K_+\sigma$,
	where

	\begin{displaymath}
	  \sigma = \begin{pmatrix} 0& 1 & 0 & 0\\ -1 &0&0&0 \\
	    0 & 0 & 1 &0\\ 0 &0 &0 &1\\
	  \end{pmatrix}
	\end{displaymath}

	Then $\K_-/\H$ is a sphere if and only if
	$\sigma\SU(3)\sigma^{-1}\cap\H=\SU(2)$(where $\SU(3)$ and $\SU(2)$ are
	the lower right blocks). We easily see that 

	\begin{displaymath}
	  \sigma\SU(3)\sigma^{-1}\cap\H=\{\diag(\bar z^2,z,\bar z^k,\bar
	  z^k)\cdot A\mid z^{2(k+1)}=1, A\in\SU(2)\}
	\end{displaymath}

	and this implies $k=0,-2$. If $k=0$, we have $N(\H)_0=\T^2\SU(2)$, which
	does not contain $\sigma$, so this is a primitive example.  If $k=-2$,
	we have $N(\H)_0=\S^1\SU(2)\SU(2)$, so up to conjugation in $N(\H)_0$ we
	have $\K_-=\K_+$, which is not possible. This finishes the case
	$\K_+=\U(3)$.
    \end{itemize}

The cases left are $(\K_+)_0=\S(\U(1)\U(1)\U(2))$ and $(\K_+)_0=\S(\U(2)\U(2))$.
Common for both consider $\SU(2)\subset\H_+$.  This would imply that $\SU(2)$
needs to act on $\K_-/\H$, and therefore $\K_-=\U(3)$ or $\K_-=\Sp(2)$, which
are cases we have considered before.

Now assume $(\K_+)_0=\S(\U(2)\U(2))$ and $\SU(2)\not\subset\H_+$. It is clear
than that $(\H_+)_0=\S^1$, and $\K_+/\H_+=\SO(4)$, implying $\SU(4)$ acts as
$\SO(6)$. Switching to $\SO(6)$, we have $(\K_+)_0=\SO(2)\SO(4)$ and
$\H_0=\SO(2)\SO(3)$, where $\SO(2)$ is the upper left and $\SO(3)$ and $\SO(4)$
are the lower right block. Since $\SO(2)$ must act on $\K_-/\H$, and the
isotropy representation of $\SO(6)/\H$ decomposes in one irreducible factor of
dimension 3,6 and 2 each, it's easy to deduce $\K_-=\SO(3)\SO(3)$, which gives a
primitive example (a Grassmanian). Note that both $l_\pm>1$, so all groups are
connected.

The last case to consider is $(\K_+)_0=\S(\U(1)\U(1)\U(2))$, where
$\SU(2)\not\subset\H_+$, as we have argued above. This implies $\K_+/\H=\mathbb
S^3$, and since $\H/\H_-$ contains $\S^1$ as a factor, we see $l_->1$, so all
groups are connected. We have $\H=\T^2$, and consider cases by the rank of
$N(\K_-)$:

\begin{itemize}
  \item Suppose $\rk(N(\K_{-}))=3$. We consider the standard representation
    $\rho$ of $\SU(4)$ on $\C^4$ and its restriction to $\K_\pm$ and $\H$. We
    know $\rho|_{\K_\pm}$ decomposes into a two-dimensional and two
    one-dimensional irreducible factors, while $\rho|_{\H}$ decomposes into four
    one-dimensional irreducible factors. 

    Consider for a moment the subcase that $\rho|_\H$ decomposes into 4
    inequivalent subrepresentations $\C e_1\oplus\C e_2\oplus\C e_3\oplus\C
    e_4$. Then the two-dimensional irreducible subspaces of $\rho|_{\K_\pm}$ are
    neccessarily given by $V_\pm=\C e_{i_\pm}\oplus\C e_{j_\pm}$ for some
    $i_\pm\ne j_\pm$. If $V_+\cap V_-\ne 0$, then after permutation we can
    assume $V_\pm\subset \C e_2\oplus\C e_3\oplus\C e_4$. But that means that
    both $\K_\pm$ are contained in the lower 3x3-block, a contradiction. But if
    $V_+\cap V_-=0$, we can use a permutation again to assume $V_+=\C
    e_1\oplus\C e_2$ and $V_-=\C e_3\oplus\C e_4$, which implies $\K_\pm\subset
    \S(\U(2)\U(2))$, again a contradiction.
    
    So we assume that $\rho|_\H$ has two equivalent one-dimensional
    subrepresentations (note that it could not be three for any embedding
    $\T^2\into\SU(4)$). This implies $N(\H)=\S(\U(2)\U(1)\U(1))$, where after a
    permutation we may assume that $\C e_1$ and $\C e_2$ are the equivalent
    subrepresentations of $\rho|_\H$. Now we can argue exactly as before, using
    that $V_\pm=A(\C e_{i_\pm}\oplus\C e_{j_\pm})$ for some $A\in N(\H)_0$. Since
    we may conjugate $\K_\pm$ by any element in $N(\H)_0$ without changing the
    manifold, we can actually assume that $V_\pm=\C e_{i_pm}\oplus\C e_{j_\pm}$
    and arrive at a contradiction as before. This finishes the subcase
    $\rk(N(\K_-))=3$ of $(\K_+)_0=\S(\U(1)\U(1)\U(2))$.

  \item The last case to consider is $\rk(N(\K_-))=2$ and
    $(\K_+)_0=\S(\U(1)\U(1)\U(2))$. This readily implies $\rk(\K_{-})=2$ and
    therefore $\K_-/\H=\mathbb S^2$. The representation $\rho|_{\K_-}$ (see the
    item before for the notation) can not be irreducible, because otherwise its
    restriction to the semisimple part would be irreducible as well, which would
    imply $\rk(\K_-)=1$. By virtue of $\rk(N(\K_{-}))=2$ it is clear that
    $\rho|_{\K_-}$ either decomposes into two 2-dimensional irreducible
    subrepresentations, or into one 3-dimensional and one 1-dimensional
    irreducible subrepresentation. This leaves us with two cases:

    \begin{itemize}
      \item[a)] \begin{displaymath}
	  \K_-=\left\{\begin{pmatrix} A & \\ & \bar A\end{pmatrix}\mid
	    A\in\U(2)\right\}
	\end{displaymath}
      \item[b)] $\K_-=\S^1\SO(3)$ where $\SO(3)$ is the upper left block and
	\begin{displaymath}
	\S^1=\diag(z,z,z,\bar z^3)
      \end{displaymath}
    \end{itemize}

    In case b), we can use  conjugation in $\SO(3)$ to achieve

    \begin{displaymath}
      \H=\{\diag(z_1,\bar z_1,1,1)\}\cdot\{\diag(z_2,z_2,z_2,\bar z_2^3)\}
    \end{displaymath}

    Thus $\rho|_\H$ decomposes into four inequivalent one-dimensional
    subrepresentations $\C e_i, i=1,\dots,4$, and the irreducible
    three-dimensional representation corresponding to $\rho|_{\K_-}$ is $\C
    e_1\oplus \C e_2\oplus\C e_3$. By primitivity, the irreducible 2-dimensional
    subspace of $\rho|_{\K_+}$ is given by $\C e_i\oplus\C e_4,i=1,\dots,3$, and
    the whole group picture is determined by $i$. Since exchanging the first two
    coordinates leaves $\K_-$ and $\H$ invariant, the cases $i=1$ and $i=2$ are
    equivalent, which leaves two possible examples. We claim that $\K_+/\H$ is
    not a sphere for $i=3$, but it is for $i=2$, which gives a primitive
    example.

    \begin{itemize}
      \item[$i=3$:] Since $\SU(2)$ acts transitively on $\K_+/\H$ with isotropy
	$\SU(2)\cap\H$, we need to show the latter is not trivial. An element
	of $\H$ being in $\SU(2)$ is equivalent to the three equations
	$z_1z_2=1$, $\bar z_1z_2=1$ and $z_2^2=1$. This is obviously true for
	$z_1=z_2=-1$, which constitutes a nontrivial element of $\SU(2)\cap \H$.
      \item[$i=2$:] Exchange the first an the third coordinate, moving $\SU(2)$
	to the lower right block and conjugating $\H$ into

	\begin{displaymath}
	  \diag(z_1,1,\bar z_1,,1)\diag(z_2,z_2,z_2,\bar z_2^3)
	\end{displaymath}

	We can read of the equations as before. This time the second coordinate
	show that any element of $\H\cap\SU(2)$ fulfills $z_2=1$, which readily
	implies $z_1=1$ from the first coordinate, so that in this case
	$\H\cap\SU(2)$ is trivial.
    \end{itemize}

    In case a), we use the fact that $-Id\in\SU(4)$ is in the kernel of the
    action, and replace $\SU(4)$ with $\SO(6)$. This gives $\K_-=\SO(2)\SO(3)$
    (where $\SO(2)$ is in the upper left and $\SO(3)$ in the lower right block)
    and $\H=\SO(2)\SO(2)$, where the two trivial subrepresentation of $\rho|_\H$
    are spanned by $e_3$ and $e_4$. By primitivity, we have $\K_+=\U(2)\SO(2)$,
    giving one primitive example, the homogeneous space

    \begin{displaymath}
      \frac{\SO(7)}{\U(2)\SO(3)}
    \end{displaymath}
\end{itemize}

\section{$\G=\SU(n),n\geq 5$}\label{sec:sun}

We claim that, up to equivalence, the diagrams of the simply connected primitive
cohomogeneity one $\SU(n)$-manifolds ($n\geq 5$) with positive Euler
characteristic are given by table \ref{table:sun}.

\renewcommand\arraystretch{1.4}

 \begin{table}[h]\caption{$\SU(n)$-cohomogeneity one manifolds for
   $n>4$}\label{table:sun}
  \begin{center}
  \begin{tabular}{|c|}\hline 
    $\S^1\SU(n-2)\subset\sigma(\U(n-1)),\U(n-1)$\\ 
    where $\S^1\SU(n-2)=\{\diag(\bar z,z,1,\dots,1)\}$ and $\sigma$ exchanges first two coordinates\\\hline
    $\S^1\SU(n-2)\subset\S(\U(2)\U(n-2)),\U(n-1)$\\
    where $\S^1=\{\diag(\bar z^{n-2},z^{n-2},\bar z,\dots,\bar z)\}\subset N(\SU(n-2))$\\\hline
    $\S^1\SU(n-2)\subset\S(\U(2)\U(n-2)),\U(n-1)$\\ 
    where $\S^1=\{\diag(\bar z^{n-2},1,z,\dots,z)\}\subset N(\SU(n-2))$\\\hline
    $\S^1\SU(n-2)\subset\SU(2)\SU(n-2),\U(n-1)$\\
    where $\S^1=\{\diag(\bar z,z,1,\dots,1)\}\subset N(\SU(n-2))$\\\hline
    $\S^1\SU(n_1-1)\SU(n_2)\subset\S(\U(n_1-1)\U(n_2+1)),\S(\U(n_1)\U(n_2))$\\
    where $\S^1=\{\diag(\underbrace{\bar z^{n_2},\dots,\bar
    z^{n_2}}_{\times n_1-1},1,\underbrace{z^{n_1-1},\dots,z^{n_1-1}}_{\times
    n_2})$\\
    and $\SU(n_1-1)\SU(n_2)\subset\H$ acts trivially on $\C n_1$
    ($n_1+n_2=n-1,n_1,n_2>1$)\\\hline
  \end{tabular}
  \end{center}
 \end{table} 
 
 The possibilities for $\K_+$ are summarized in table \ref{table:sunk+} (again,
 we refer to section \ref{sec:procedure}).

 \begin{table}[h]\caption{Possibilities for $K_+$}\label{table:sunk+}
  \begin{center}
  \begin{tabular}{|c|c|l|l|}\hline 
    Factors& Subcase & Group &Conditions\\\hline 
    2 & 2a &$\U(n-1)$ & $n>2$ \\\hline 
    3 & 3a &$\S(\U(n_1)\U(n_2))$ & $n_1, n_2 >1, n_1+n_2=n$\\\hline 
    3 & 3b &$\S^1\U(n-2)$ & $n>3$ \\\hline
    4 & 4a &$\S(\U(1)\U(n_1)\U(n_2))$ & $n_1, n_2>1, n_1+n_2=n-1$\\\hline 
  \end{tabular} 
  \end{center}
 \end{table} 
    
\begin{itemize} 
  
  \item[2a)] In this case we assume $\K_+=\U(n-1)$, where we can assume
    $\SU(n-1)$ is the lower right block (this also determines the center of
    $\U(n-1)$). If $\SU(n-1)\subset \H_+$, we have
    $\SU(n)\subset\K_-$, which is treated in \ref{subsubsec:fixedpoint}. If
    $n>5$ and $\SU(n-1)\not\subset\H_+$, we can use table
    \ref{table:sphereactions} to see that $\K_+/\H_+$ is $\U(n-1)$, so by
    \ref{subsubsec:un1spheres} we have $\H=\S^1_k\SU(n-2)$ where
    
    \begin{displaymath}
      \S^1_k=\{\diag(\bar z^{n-2},z^{(k+1)(n-2)},z^{-k},\dots,z^{-k})\mid z\in
      \S^1\}
    \end{displaymath}

    For $n=5$, there is the additional possibility that $\SU(4)=\Spin(6)$ acts
    with isotropy $\Spin(5)=\Sp(2)$, i.e.\ $\H=\S^1\Sp(2)$. But
    $\K_+/\H_+=\SO(6)$ implies $\S^1\subset\H_+$, and $\K_+/\S^1=\SO(6)/\{\pm
    Id\}$, which by table \ref{table:sphereactions} cannot act transitively on
    any sphere. So the above claim holds for $n=5$ as well. We differentiate
    between the possible values for $\H_-$ (the kernel of the action of $\K_-$
    on $\K_-/\H$).

    If $\K_-$ acts almost effectively on $\K_-/\H$, we have $\K_-\simeq\U(n-1)$, and
    by studying the isotropy representation of $\H$ in $\G$ we see that
    primitivity implies $K=\sigma K_+\sigma^{-1}$ where     

    \begin{displaymath}
      \sigma=
      \begin{pmatrix}
	0& 1 & \\
	-1& 0 & \\
	 &  & 1 \\
	 & & & \ddots \\
	 &&&&1
      \end{pmatrix}
    \end{displaymath}
    
    where the empty spaces are filled up with zeroes. $\H$ can only occur as an
    isotropy group for a transitive almost effective action on a sphere of
    $\K_-$ in the case $k=0$ or $k=-2$. In the latter case, we have
    $N(\H)_0=\SU(2)\H$, which contains $\sigma$, so this is not a primitive
    example. In the former case, $\K_-$ and $\K_+$ are not conjugate by an
    element of $N(\H)_0=\T^2\SU(n-2)$, so we get a primitive example.

    The next case is $\H_-=\S^1_k$. Either $\S^1_k$ intersects $\SU(n-2)$ and
    $\H/\H_-$ is a proper quotient of $\SU(n-2)$, and $\K_-/\H_-=\SO(7)$, which
    follows from studying the classification of transitive effective actions on
    spheres given in table \ref{table:sphereactions}. Since $N(\S^1_k)$ does
    not contain $\Spin(7)$ or $\SO(7)$, this is impossible. The other
    possibility is that $\S^1_k$ does not intersect $\SU(n-2)$, which is the
    case exactly if $n-2$ divides $k$. Since $\SU(n-1)\subset N(\S^1_k)$, we
    have $k=n-2$, which is impossible because $\SU(n-1)\cap\S^1_{n-2}\ne 0$.
    
    If $(\H_-)_0=\SU(n-2)$, we have $\H/\H_-\simeq \S^1$ and
    $\K_-/\H\in\{\mathbb S^2,\mathbb S^3\}$. Also $\K_-\subset N(\SU(n-2)) =
    \S^1\SU(2)\SU(n-2)$. If $\K_-=\SU(2)\SU(n-2)$, that implies $k=0$ and gives
    a primitive example ($\SU(n)$ acting on $\mathbb HP^{n-1}$). If
    $\K_-=N(\SU(n-2))$, the condition that $\SU(2)\cap \H$ is trivial is

    \begin{displaymath}
      \SU(2)\cap\H=\{\diag(\bar z^{n-2},z^{n-2},1,\dots,1)\mid
      (z^k)^{n-2}=1\}\text{ is trivial }
    \end{displaymath}

    which implies $|k|=1$ and gives 2 primitive examples ($k=-1$ gives the
    Grassmanian $\SU(n+1)/\S(\U(2)\U(n-1))$).
 \item[3a)] For the case $\K_+=\S^1\SU(n_1)\SU(n_2)$, we will first do the most
   general case and after that care for the exceptional actions. Since
   $n_1,n_2>1$, one of $\SU(n_1), \SU(n_2)$ must act on $\mathbb S^{l_+}$, and
   we may assume it's $\SU(n_1)$.

   Let's first show that without any further assumptions, we  have
   $\Sp(2)\not\subset \H$. Otherwise $n_1=4$ (and
   $\S^1\SU(4)\ne\U(4)$ as well, but we don't need that), and
   $\H=\S^1\Sp(2)\SU(n-4)$. This would mean $\S^1\SU(n-4)\subset\H_-$, but
   $N(\Sp(2))_0=\H$, so $\K_-/\H$ could not possibly be a sphere of positive
   dimension.

   Now assume $n_1>3,n_2>2$. By the previous paragraph the action of $\K_+$ on
   $\K_+/\H$ is given as follows: The matrix $\begin{pmatrix}A & 0\\ 0&
     B\end{pmatrix}$ with $A\in \U(n_1),B\in\U(n_2), \det A\det B=1$ acts as
     $(\det A)^kA$, giving

    \begin{displaymath}
      H=\left\{\begin{pmatrix} A' &  & \\ & a & \\ &    & B\\
    \end{pmatrix}\mid a\det A'\det B=1 \text{ and } a^{k+1}=\overline{\det
    A'}^k\right\} 
    \end{displaymath}
  
    for some $k\in\Z$.  Now $\SU(n_2)$ cannot act trivially on $\K_-/\H$. If
    $n_2=4$, we could have $\K_-/\H_-\simeq\SO(7)$, but there is no $\SO(7)$ or
    $\Spin(7)$ in the normalizer $N(\H)=\S(\U(n_1-1)\U(n_2))$. Since $n_2>2$, by
    the classification of effective actions on spheres we have
    $K_-=\S^1\SU(n_1-1)\SU(n_2+1)$, where $\SU(n_1-1)$ is the upper left and
    $\SU(n_2+1)$ the lower right block. The action of $\K_-$ on $\K_-/\H$ is given
    in a similar fashion as before, giving the isotropy group 
  
  \begin{displaymath} \tilde H=\left\{\begin{pmatrix} A'&  & \\ & a & \\ &    &
    B\\ \end{pmatrix}\mid a\det A'\det B=1 \text{ and } a^{l+1}=\overline{\det
    B}^l\right\} \end{displaymath}
    
    for some $l\in \Z$. We have $\tilde \H=\H$ only in the case $k=l=0$ as
    follows: For any matrix in $\H\cap \tilde \H$ we have

    \begin{align*} \det \bar A'&=a \det B\Rightarrow a^{k+1}=a^k\det
      B^k\Rightarrow a=\det B^k\\ \det \bar B &=a \det A'\Rightarrow
      a^{l+1}=a^l\det A'^l\\ &\Rightarrow \det B^{k(l+1)}=\det \bar B^l
      \Rightarrow \det B^{k(l+1)+l}=1\\ &\Rightarrow k(l+1)+l=0 \Rightarrow
      k=l=0 \text{ because } k,l\in\Z \end{align*}

    giving the only primitive example in this case, the Grassmanian
    $\SU(n_1+n_2+1)/\S(\U(n_1)\U(n_2+1))$.

    Now assume $n_1=2$, which by $n\geq 5$ implies $n_2\geq 3$. Since
    $\text{corank}(\H)=1$ (see lemma \ref{cor:corH}), we have $l_+=3$. If
    $\K_+/\H=\mathbb S^3$, we have $\H=\S^1\SU(n_2)$. If $n_2=4$, we could again
    have $\K_-/\H_-=\SO(7)$, but there's no $\SO(7)$ or $\Spin(7)$ in $\SU(6)$.
    This shows $\K_-=\U(n-1)$, which is a case we treated before.

    The next case is $n_1=3$ and $n_2\geq 2$. But $n_2>2$ actually implies that
    the argument given in the beginning of this section applies as well, so we
    only need to deal with the case $n_2=2$. Then $\H=\S^1_k\SU(2)\SU(2)$, where

    \begin{displaymath}
      \S^1_k=\{\diag(z^{2-2k},z^k,z^k,\bar z,\bar z)\mid z\in\S^1\}
    \end{displaymath}

    If $\K_-/\H=\mathbb S^5$, then again the argument above yields one primitive
    example (corresponding to $k=1$ here). From the classification of transitive
    actions on a sphere and the fact that $\SU(2)\subset \H/\H_-$, we are left
    with 2 further possibilities: $\K_-/\H_-=\SO(5)$ and $\K_-/\H_-=\SO(4)$.In
    the first case we can deduce $\K_{-}=\S^1\Sp(2)$, which implies $k=-1$ and
    results in $\mathbb CP^9$. In the latter case, we would have
    $\SU(2)\subset\H_-$, but there is no $\SO(4)$ or
    $\Spin(4)=\SU(2)\times\SU(2)$ in the normalizer $N(\SU(2))$ in $\SU(5)$.

    The last case to consider is $n_1>3, n_2=2$. We will argue that there is no
    additional possibility to the one given at the beginning of this section. We
    have $\H=\S^1\SU(n-3)\SU(2)$, and since $n-3\geq 3$ we have
    $\SU(n-3)\subset\H_-$. But $N(\SU(n-3))=\S^1\SU(3)$, so the only possibility
    is $\SU(3)\subset\K_-/\H_-$, which is the aforementioned argument.

      \item[4a)] We assume $\SU(n_1)$ is the upper left, $\SU(n_2)$ the lower
	right block. Note that this determins the center $\T^2$ of $\K_+$. By
	remark \ref{rem:numberoffactors} and the fact that in this case $\K_+$
	has 4 factors, we can assume $(\H_{-})_0=\S^1\SU(n_2)$ and therefore
	$\H=\T^2\SU(n_1-1)\SU(n_2)$, where $\SU(n_1-1)$ is the upper left and
	$\SU(n_2)$ is the lower right block. If we assume
	$\K_-\simeq\U(n_1-1)\U(n_2+1)$ and $n_1>2$ then we know the resulting
	manifold will not be primitive for the following reason:
	$\K_-\into\SU(n)$ induces a 1-dimensional subrepresentation, which by
	the choice of $\K_+$ and $\H$ is necessarily given by $e_{n_1}$ or
	$e_{n_1+1}$. The first case implies
	$\K_\pm\subset\S^1\SU(n_1)\SU(n_2+1)$, and the second
	$\K_\pm\subset\U(n-1)$.

	We divide the remainig cases:

	\begin{itemize}
	  \item[$\bullet$] If $n_1=2$, then $\H=\T^2\SU(n_2)$, and
	    $\H\into\SU(n)$ induces 3 one-dimensional subrepresentations via the
	    standard representation. We know $\K_-\simeq\T^2\SU(n_2+1)$, for
	    which we only need to note that $\T^2\Sp(2)$ is not a subgroup of
	    $\SU(4)$ by table \ref{table:maxranksubgroups}. Now
	    $\K_-\into\SU(n)$ induces 2 one-dimensional subrepresentations, and
	    if one of those is given by $e_3$, we have $\K_\pm\subset\U(n-1)$.
	    But in the other case $\K_\pm\subset\S^1\SU(2)\SU(n-2)$, so there is
	    no primitive example in this case.
	  \item[$\bullet$] We will now argue that indeed
	    $\K_-\simeq\U(n_1-1)\U(n_2+1)$ in all other cases. By primitivity we
	    know $\left(\H/\H_-\right)_0=\S^1\SU(n_2)$, so by the classification
	    given in table \ref{table:sphereactions} this is true for $n_1\geq
	    3,n_2\geq 3$. Having already dealt with the case $n_1=2$, we only
	    need to consider $n_1\geq 3,n_2=2$, but again we only need to note
	    that $\T^2\Sp(2)\SU(n_1-4)$ is not a subgroup of $\SU(n)$ by table
	    \ref{table:maxranksubgroups}. This shows
	    $\K_-\simeq\U(n_1-1)\U(n_2+1)$ as desired, and we know that the
	    manifold is not primitive by what was said above.  
	\end{itemize}

  \item[3b)] If $\SU(n-2)\subset\H_+$, we have $\SU(n-1)\subset\K_-$ (remember
    that $n\geq 5$), so that $\K_-\simeq \U(n-1)$, which was treated before. So
    now assume $\SU(n-2)\not\subset\H_+$, which implies $\H=\T^2\SU(n-3)$. We
    distinguish the cases $\SU(n-3)\not\subset\H_-$ and its opposite.

    \begin{itemize}
      \item[$\bullet$] If $\SU(n-3)\not\subset\H_-$, we have $\K_-\simeq\K_+$.
	The factor $\SU(n-2)$ of $\K_-$ then has a 2-dimensional trivial
	subrepresentation when consideres as a subgroup of $\SU(n)$ which is
	necessarily given by
	$\C e_i\oplus \C e_j$ for $i,j\in\{1,2,3\}$:
	\begin{itemize}
	  \item If $i=1,j=2$, we have $\K_-=\K_+$, so the example is not
	    primitive.
	  \item If $i=1,j=3$, we have $\K_\pm\subset\S(\U(1)\U(n-1))$, so again
	    this is not primitive
	  \item If $i=2,j=3$, exchange the coordinates $e_1$ and $e_2$, after
	    which we are in the previous case again. Note that this exchange
	    might change every group of the diagram, yet still shows that it is
	    not primitive.
	\end{itemize}
      \item[$\bullet$] If $\SU(n-3)\subset\H_-$, we have $\H/\H_-\simeq \S^1$,
	which implies $\K_-/\H_-$ is one of $\U(2), \SO(3)$. In the first case
	$\K_-\simeq\S(\U(2)\U(n-3))$, which was previously treated in case 4a.
	In the latter case, we have $\K_-\simeq\SO(3)\U(n-3)$, and we conjugate
	the whole diagram to make $\K_-$ standard: $\SO(3)$ is the upper left
	block, $\SU(n-3)$ the lower right and $\S^1$ is diagonally embedded, so
	it commutes with both $\SO(3)$ and $\SU(n-3)$ but does not act trivially
	on any of the $\C e_i$ for $i\in\{1\dots,n\}$. We can assume $\H$ is
	then given by $\S^1_1\S^1_2\SU(n-3)$ where

	\begin{displaymath}
	  \S^1_1=\{\diag(z_1,1,\bar z_1,1,\dots,1)\mid z_1\in\S^1\}
	\end{displaymath}

	and

	\begin{displaymath}
	  \S^1_2=\{\diag(\bar z_2^{n-3},\bar z_2^{n-3},\bar
	  z_2^{n-3},z_2^3,\dots,z_2^3)\mid z_2\in\S^1\}
	\end{displaymath}

	Now $\K_+=\S(\U(1)\U(1)\U(n-2))$ is determined by $\SU(n-2)$, which
	again has a 2-dimensional trivial subrepresentation when restricting the
	standard representation of $\SU(n)$, which is necessarily given by $\C
	e_i\oplus\C e_j$ for $i,j\in\{1,2,3\}$ and $i\ne j$. Note that the cases
	$i=1,j=2$ and $i=2,j=3$ are equivalent, and we claim this gives a
	(primitive) example. The last case $i=1,j=3$ does not give an example
	for $\K_+/\H$ is not a sphere:

	\begin{itemize}
	  \item If $i=1,j=2$, an element of $\H$ that is also in $\SU(n-2)$
	    satisfies $z_1\bar z_2^{n-3}=1$ and $\bar z_2^{n-3}$ from the first
	    two coordinates, which implies $z_1=1$ and so
	    $\diag(1,1,1,z_2^3,\dots,z_2^3)\in\SU(n-3)$, which shows
	    $\SU(n-2)\cap\H=\SU(n-3)$, so $\K_+/\H=\SU(n-2)/(\H\cap\SU(n-2))$ is
	    a sphere.
	  \item If $i=1,j=3$, an element of $\H$ is in $\SU(n-2)$ if and only if
	    $z_1\bar z_2^{n-3}=1$ and $\bar z_1\bar z_2^{n-3}$ from the first
	    and the third coordinate. Choose any $z_2$ such that $z_2^{n-3}=-1$
	    and $z_1=-1$, so that these equations are fulfilled. But then note
	    that the second coordinate of this element is $-1$, which shows that
	    it is in $\H\cap\SU(n-2)$, but not in $\SU(n-3)$. Therefore
	    $\K_+/\H=\SU(n-2)/(\H\cap\SU(n-2))$ is not a sphere.
	\end{itemize}
    \end{itemize}

\end{itemize}
\section{$\G=\SO(2n+1)$, $n\geq 3$}\label{sec:so2n+1}

Up to equivalence, the simply connected primitive cohomogeneity one
$\SO(2n+1)$-manifolds ($n\geq 3$) are given by table \ref{table:so2n+1}.

 \begin{table}[h!]\caption{$\SO(2n+1)$-cohomogeneity one manifolds for
   $n>2$}\label{table:so2n+1}
  \begin{center}
  \begin{tabular}{|c|}\hline 
    $\SO(2n-1)\subset\SO(2)\SO(2n-1),\SO(2n)$\\\hline
    $\O(2n-1)\subset\SO(2)\SO(2n-1),\O(2n)$\\\hline
    $\SU(3)\subset\G_2,\U(3), n=3$\\\hline
    $\S^1\SU(3)\subset\S^1\G_2,\U(4), n=4$\\\hline
    $\SO(2n_1+1)\SO(2n_2-1)\subset\SO(2n_1+2)\SO(2n_2-1),\SO(2n_1+1)\SO(2n_2)$\\
    where $n_1+n_2=n$\\\hline
    $\SO(2)\SO(2n-3)\subset\U(2)\SO(2n-3),\SO(2)\SO(2n-2)$\\\hline
    $\SO(2n_1-1)\U(n_2)\subset\SO(2n_1-1)\U(n_2+1),\SO(2n_1)\U(n_2)$\\
    where $n_1+n_2=n$\\\hline
    $\T^2\SU(n-2)\subset\SO(3)\S^1\SU(n-2),\SO(2)\U(n-1)$\\
    where $\S^1=\{\diag(1,1,1,z,\dots,z)\},
    \T^2=\{\diag(1,z,z^2,1,\dots,1)\}\cdot \S^1$\\
    and $\SO(3)\into\SO(5)$ is irreducible ($\SO(5)$ being the upper left
    block)\\\hline
  \end{tabular}
  \end{center}
 \end{table} 
 
The possibilities for $\K_+$ are given by table \ref{table:so2n+1k+}. 

\begin{table}[h]\caption{Possibilities for $K_+$}\label{table:so2n+1k+}
  \begin{center}
    \begin{tabular}{|c|c|l|l|}\hline
      Factors& Subcase &Group &Conditions\\\hline
      1 & 1a &$\SO(2n)$ & -- \\\hline
      2 & 2a &$\SO(2n_1)\SO(2n_2)$ & $n_1, n_2 \geq 3, n_1+n_2=n$\\\hline
      2 & 2b &$\SO(2n_1+1)\SO(2n_2)$ & $ n_2\geq 2, n_2+n_2=n$\\\hline
      2 & 2c &$\SO(2)\SO(2n-2)$ & $n\geq 4$\\\hline
      2 & 2d &$\SO(2)\SO(2n-1)$ & --\\\hline
      2 & 2e &$\U(n)$ & -- \\\hline
      3 & 3a &$\SO(2n_1)\U(n_2)$ & $n_1, n_2\geq 2, n_1+n_2=n$ \\\hline
      3 & 3b &$\SO(2n_1+1)\U(n_2)$ & $n_2\geq 2, n_1+n_2=n$\\\hline
      3 & 3c &$\SO(2)\U(n-1)$ & -- \\\hline
      3 & 3d &$\SO(2n-4)\SO(4)$ & $n\geq 5$\\\hline
      3 & 3e &$\SO(2)\SO(4)$ & $n=3$\\\hline
      3 & 3f &$\SO(2n-3)\SO(4)$ & -- \\\hline
      4 & 4a &$\U(n_1)\U(n_2)$ & $n_1,n_2\geq 2, n_1+n_2=n$\\\hline
      4 & 4b &$\U(n-2)\SO(4)$ & $n \geq 4$\\\hline
      4 & 4c &$\SO(4)\SO(4)$ & $n=4$\\\hline
    \end{tabular}
  \end{center}
\end{table}

If a group of complex matrices is involved (e.g. $\U(n)$), we will deliberatly
use complex notation for the corresponding real matrices. In particular, for
$e^{i\varphi}=z\in S^1$, we will use $\diag(z,\ldots,z)$ for the matrix
containing $\begin{pmatrix}\cos\varphi &\sin\varphi \\ -\sin\varphi
&\cos\varphi\end{pmatrix}$ on the diagonal $2\times2$-blocks and 0 everywhere
else.

\begin{itemize}

  \item[1a)] We have $\H=\SO(2n-1)$. There are two possibilities for $\K_-$,
    namely $\K_-\simeq \SO(2n)$ and $\K_-\simeq\SO(2)\SO(2n-1)$. The first
    choice does not lead to a primitive manifold, since 
    $\K_-$ is conjugate to $\K_+$ via the matrix

    \begin{displaymath}
      \begin{pmatrix}
	\cos \varphi & \sin\varphi &  &  \\
	-\sin\varphi & \cos\varphi &  &\\
	 & & & & \\
	 & & & E_{2n-1}
      \end{pmatrix}
    \end{displaymath}
  
    which is in $N(\H)_0$. There is only one $\SO(2)$ in the normalizer of
    $\SO(2n-1)$, so the second choice leads to exactly one primitive manifold, a
    Grassmanian.

    We do have $l_-=0$ here, so $\K_+$ might be non-connected. It is clear then
    that $\H=\Z_2\SO(2n-1)\simeq \O(2n-1)$ and $\K_+=\Z_2\SO(2n)\simeq\O(2n)$.
    This is $\SO(2n+1)$ acting on $\mathbb CP^{2n}$.
  \item[2a)] We can assume $\H=\SO(2n_1-1)\SO(2n_2)$, and it is clear that
    $\K_-\simeq\SO(2n_1-1)\SO(2n_2+1)$. There are two possibilities for
    $\SO(2n_2+1)$ in the normalizer of $\SO(2n_1-1)$, both of which are
    conjugate by a change of coordinate, which can be achieved by conjugation
    with a matrix similar to the one given in 1a, which is in
    $N(\H)_0=\SO(2)\SO(2n_1-1)\SO(2n_2)$. But it is obvious for at least one of
    the two possibilities that $\K_\pm\subset\SO(2n)$.

  \item[2b)] Since $\text{corank}(\H)=1$ by corollary \ref{cor:corH}, we have
    $\H=\SO(2n_1-1)\SO(2n_2+1)$, implying  $\K_-=\SO(2n_1-1)\SO(2n_2+2)$. This
    gives one primitive example, a Grassmanian.
  \item[2c)] First assume $\H=\SO(2)\SO(2n-3)$. Further assuming
    $\K_-/\H_-=\U(2)$, we see $\K_-\simeq\U(2)\SO(2n-3)$. The choices for
    $\U(2)$ in the normalizer of $\SO(2n-3)$ are given by the center, which is
    given by either $\{\diag(z,z,1,\dots,1)\mid z\in\S^1\}$ or $\{\diag(\bar
    z,z,1,\dots,1)\mid z\in\S^1\}$. But complex conjugation of the first
    component in $\T^2\subset\SO(4)$ is given by conjugation with
    $\diag(-1,1,\dots,1)$, which leaves $\H$ invariant as well as $\K_+$, so we
    only get one new example. It is the homogeneous space

    \begin{displaymath}
      \frac{\SO(2n+2)}{\SO(2n-2)\U(2)}
    \end{displaymath}

    The next possibility is $\K_-\simeq\SO(3)\SO(2n-3)$, and again from checking
    the isotropy representation there are 2 choices for $\SO(3)$, but both are
    conjugate by a change of two coordinates, which is in $N(\H)_0$ as in 1a.
    The result is not primitive.

    If $\H=\Z_k\SO(2n-2)$, we have $\K_-=\Z_k\SO(2n-1)$. In the case that $\H$
    is connected, there are several choices for $\K_-=\SO(2n-1)$, but all of
    them are conjugate in $N(\H)_0=\SO(3)\SO(2n-2)$, and obviously not
    primitive. If $\H$ is not connected, there's only one choice for
    $\SO(2n-1)$, but again $\K_\pm\subset\SO(2)\SO(2n-1)$.

  \item[2d)] We only need to note $\SO(2n-1)\subset\H_+$, since
    $\text{corank}(\H)=1$ by corollary \ref{cor:corH}, which leads to
    $(\K_-)_0=\SO(2n)$ which was treated in 1a.
  \item[2e)] If $\SU(n)\subset\H_+$, then $n=3$ or $n=4$, because for $n\geq 5$
    we would have $\SU(n+1)\subset \K_-$, but there is no embedding of
    $\SU(n+1)\into \SO(2n+1)$ (as seen easily from checking representations). If
    $n=3$, there is the possibility that $\H=\SU(3)$ and $\K_-=\G_2$, giving one
    primitive example: It is primitive because $\K_+$ is maximal and does not
    contain $\K_-$ or any of its conjugates, and we can see there is only one in
    the following way. We conjugate the diagram so that $\K_-$ is a given, fixed
    subgroup of type $\G_2$ that contains the standard $\SU(3)$ lower right
    block. This already determines $\K_+$ as the centralizer of $\H$. Now if
    $n=4$, we can have $\H=\SU(4)$ and $\K_-=\Spin(7)$, which is not primitive
    because $\K_\pm\subset\SO(8)$.
    
    From now on we assume $\H=\S^1_k\SU(n-1)$, where

    \begin{displaymath}
      \S^1_k=\{\diag(1,z^{(k+1)(n-1)},z^{-k},\dots,z^{-k})\mid z\in\S^1\}
    \end{displaymath}

    and $\SU(n-1)$ is in the lower right block.

    First assume $k=-1$. We have $N(\H)_0=\SO(3)\SU(n-1)$. If $\K_-\simeq\U(n)$,
    there are several choices, but all of them can be conjugated into $\SO(2n)$
    by a simple change of coordinates, which is in $N(\H)_0$, so there's no
    primitive example. We cannot have $\SU(n-1)\subset\H_-$, but
    $\S^1\not\subset\H_-$, since there are no subgroups in $N(\SU(n-1))$ of type
    $\U(2), \SO(3), \SU(2)$ containing $\S^1_k$. Lastly $\U(n-1)\subset\H_-$
    implies $\K_-\simeq\SO(2)\U(n-1)$, and again conjugation in $N(\H)_0$ leads
    to $\K_-\subset\K_+$.  

    Now for $k\ne-1$ we have $N(\H)_0=\SO(2)\U(n-1)$. If $\K_-\simeq\U(n)$ there
    are two choices for embedding $\U(n)$, but both are in $\SO(2n)$ as well as
    $\K_+$. We cannot have $\H=\H_-$, because $N(\H)_0\subset\U(n)$, so we
    would not get a primitive manifold. For the same reasons as before,
    $\SU(n-1)\subset \H_-$, but $\S^1_k\not\subset\H_-$ is not possible. If
    $\H_-=\S^1$, we have $n=4$ and $\K_-=\S^1\G_2$ or $n=5$ and
    $\K_-=\S^1\Spin(7)$, for there is no embedding
    $\SU(n+1)\into\SO(2n+1)$. In the latter case $\K_\pm\subset\SO(10)$, so the
    result is not primitive. The first case results in one primitive example:
    $\K_+$ is maximal and does not contain any conjugate of $\K_-$, so it's
    primitive, and we can see it's the only example in the same way as the
    beginning of the section: Conjugating $\K_-,\H$ in a standard form
    determines $\K_+$ uniquely.

  \item[3a)] If $\H=\SO(2n_1-1)\U(n_2)$, it is clear that
    $\K_-\simeq\SO(2n_1-1)\U(n_2+1)$, and there are 2 choices for embedding
    $\K_-$ into $\SO(2n+1)$ while containing $\H$, corresponding to the
    embedding of the additional coordinate of its center over the one of $\H$.
    Both are conjugate by 

    \begin{displaymath}
      \diag(\underbrace{1,\dots,1}_{2n_1-1},-1,\underbrace{1,\dots,1}_{2n_2+2})
    \end{displaymath}

    where we assume that $\U(n_2)$ is embedded as the lower right block, and
    $\SO(2n_1-1)$ in the upper left. This conjugation is not in
    $N(\H)_0=\SO(2)\H$, but restricted to $\K_-$ it is the same as conjugation
    with
    
    \begin{displaymath}
      \diag(\underbrace{-1,\dots,-1}_{2n_1-1},-1,\underbrace{1,\dots,1}_{2n_2+2})
    \end{displaymath}

    which is in $\SO(2n_1)$ and leaves both $\K_+$ and $\H$ invariant. So we
    only obtain one example, which is primitive, because $\K_-$ is a maximal
    subgroup of maximal rank not containing and not contained in $\K_+$ (and not
    isomorphic to it). This is the homogeneous space

    \begin{displaymath}
      \frac{\SO(2n+2)}{\SO(2n_1)\U(n_2+1)}
    \end{displaymath}

    The second possibility is $\H=\SO(2n_1)\U(n_2-1)_k$ (where
    $\U(n_2-1)_k=\S^1_k\SU(n_2-1)$ similar to 2a). It is clear that
    $\K_-\simeq\SO(2n_1+1)\U(n_2-1)_k$, and for $k\ne -1$ there is just one
    possibility for that, which is not primitive, for $\K_\pm\subset
    \SO(2n_1+1)\SO(2n_2)$. If $k=-1$, there are actually 3 possibilities, but
    all again differ only by a change of coordinates, which can be done in
    $N(\H)_0=\SO(3)\H$, so again no primitive example arises.

 \item[3b)] We cannot have $\H_+=\U(n_2)$, for that would imply
   $\SO(2n_1)\U(n_2+1)\subset\K_-$, but this can't be embedded into $\SO(2n+1)$.
   Also, $\H_+\ne\Z_k\SU(n_2)$ for any $k$, because $\SO(2n_1+1)\SU(n_2)$ is not
   isotropy group of any almost effective transitive action on a sphere.

   So we have $\H=\SO(2n_1+1)\U(n_2-1)_k$ (see 3b for notation), and therefore
   $\SO(2n_1+2)\subset\K_{-}$ is in the normalizer of $\U(n_2-1)_k$, which is
   only possible for $k=-1$. If $k=-1$, we have 2 possibilities for $\K_-$, both
   conjugate by a change of coordinates, which is in $\N(H)_0=\SO(2)\H$. We are
   left with $\K_-=\SO(2n_1+2)\U(n_2-1)$, which is primitive as shown in 3a.

 \item[3c)] First suppose $\SU(n-1)\subset\H_+$. This implies
   $\SU(n)\subset\K_-$ by primitivity and the classification of transitive
   actions on spheres. Therefore $(\K_{-})_{0}\simeq\U(n)$, which is contained
   in case 2e.  From  now on we assume $\SU(n-1)\not\subset\H_+$.

   By remark \ref{rem:numberoffactors} we can deduce $\H_0=\T^2\SU(n-2)$.
   Suppose $\SU(n-2)\not\subset\H_-$ (note that this implies $n>3$). Then
   $\K_-/\H_-\simeq \S^1\SU(n-1)$. Now we can deduce
   $\K_\pm\subset\SO(2n)$ as follows: All irreducible real representations of
   $\K_\pm$ and $\H$ are even-dimensional, so there is at least one
   1-dimensional trivial representation given by the embedding
   $\K_-\into\SO(2n+1)$, which of course stays trivial when restricted to $\H$.
   If the embedding $\H\into\SO(2n+1)$ induces only one such representation,
   then this is neccessarily the same as the one for $\K_+$, which shows the
   claim. If $\H\into\SO(2n+1)$ induces three trivial 1-dimensional
   representations, then $N(\H)_0\supset\SO(3)$, and we may again assume
   $\K_-\subset\SO(2n)$ by section \ref{subsec:uniq}.

   Now we are left with the case $\H=\T^2\SU(n-2)$ and $\SU(n-2)\subset\H_-$,
   which implies $l_-=2$ or $l_-=3$. We divide cases by the possibilities for
   $\K_-$:

   \begin{itemize}
     \item[$\bullet$] Suppose $l_-=2$, $\K_-\simeq\SO(3)\S^1\SU(n-2)$ and
       $\SO(3)\into\SO(2n+1)$ is irreducible in $\SO(5)\subset N(\S^1\SU(n-2))$. We
       note that $\SO(2)\subset\SO(3)$ has weights 1 and 2 in $\SO(5)$ in
       this case, and therefore we have $\H=\S^1_1\S^1_2\SU(n-2)$, where
       $\SU(n-2)$ is the lower right block and 

       \begin{eqnarray*}
	 \S^1_1&=\{\diag(1,z,z^2,1,\dots,1)\mid z\in\S^1\}\\
	 \S^1_2&=\{\diag(\underbrace{1,\dots,1}_{5},\underbrace{z,\dots,z}_{n-2}\mid z\in\S^1\}
       \end{eqnarray*}

       This yields one primitive example.

     \item[$\bullet$] If $l_-=2$, $\K_-\simeq\SO(3)\S^1\SU(n-2)$ and
       $\SO(3)\into\SO(5)\subset N(\SU(n-2))$ has a trivial 2-dimensional
       subrepresentation, then we can reconjugate the diagram to achieve
       $\H=\S^1_1\S^1_2\SU(n-2)$ where $\SU(n-2)$ is the lower right block and

       \begin{eqnarray*}
	 \S^1_1&=\{\diag(1,z,1,\dots,1)\mid z\in\S^1\}\\
	 \S^1_2&=\{\diag(1,1,1,z^{l_2},z^{k_2},\dots,z^{k_2})\mid z\in\S^1\}
       \end{eqnarray*}
       where $l_2,k_2\in\Z$. The trivial 1-dimensional subrepresentation of
       $\K_+\into\SO(2n+1)$ is then neccessarily given by $\R e_1$, which shows
       $\K_\pm\subset\SO(3)\U(n-2)$, and there is no primitive example in this
       case.

     \item[$\bullet$] If $l_-=2$ and $\K_-\simeq\SU(2)\S^1\SU(n-2)$, we have
       $\K_\pm\subset\SO(2n)$ by the same reasoning that we used in the first
       paragraph of this case.

     \item[$\bullet$] If $l_-=3$, then $\K_-\simeq\U(2)\U(n-2)$ and again
       $\K_\pm\subset\SO(2n)$ by the same reasoning that we used in the first
       paragraph of this case.
   \end{itemize}
 \item[3d-f)] All these are subject to the considerations in 2a and 2b, since
   the corank of $\H$ in $\G$ is 1 by corollary \ref{cor:corH}, so that
   $\K_+/\H$ is not $\S^2$.
 \item[4a)] We can assume $\H=\T^2\SU(n_1-1)\SU(n_2)$ and therefore
   $\K_-\simeq\T^2\SU(n_1-1)\SU(n_2+1)$. We will identify $\K_-$ via its center. It
   has to commute with $\SU(n_1-1)\SU(n_2)$, so we need to look at embeddings
   $\S^1\into\SO(3)$, all of which are conjugate in $\SO(3)$ up to complex
   conjugation. But this implies that $\K_\pm\subset\SO(2n)$ after conjugating
   $\K_-$ with an element of $N(\H)_0$, which contains the $\SO(3)$ in question.

 \item[4a-b)] Again, we have $\text{corank}(\H)=1$ by corollary \ref{cor:corH},
 so in particular $\K_+/\H\ne\S^2$. All the other cases have been treated
 before.\end{itemize}

\subsection{$\G=\Spin(2n+1)$, $n\geq 3$}\label{subsec:spin2n+1}

Up to equivalence, the simply connected primitive cohomogeneity one
$\Spin(2n+1)$-manifolds ($n\geq 3$) are given by table \ref{table:spin2n+1}.

\begin{table}[h]\caption{$\Spin(2n+1)$-cohomogeneity one manifolds for
  $n>2$}\label{table:spin2n+1}
  \begin{center}
  \begin{tabular}{|c|}\hline 
	$\S^1_k\SU(2)\subset\U(3),\SO(2)\SO(5), n=3$\\
	where $\S^1_k=\{\diag(z^2,z^k,z^k,1)\}$, $k=1,-3$ and\\
	$\SU(2)$ acts trivially on $\R e_1,\R e_2,\R e_7$\\\hline
	$\SU(3)\subset\G_2,\SO(6), n =3$\\\hline
	$\SU(3)\subset\G_2,\U(3), n=3$\\\hline
	$\SU(2)\S^1_\H\SU(2)\subset\U(2)\SO(5),\SO(5)\U(2),n=4$\\
	where $\S^1_\H=\{\diag(z^{l_1},z^{l_1},1,z^{l_2},z^{l_2})\}$ and
	$_1=2,l_2=1$ or $l_1=1=l_2$.\\\hline
  \end{tabular}
  \end{center}
 \end{table} 

The possibilities for $\K_+$ are given by table \ref{table:spin2n+1k+}.

\begin{table}[h]\caption{Possibilities for $K_+$}
  \label{table:spin2n+1k+}
  \centering
    \begin{tabular}{|c|c|l|l|}\hline
      Factors& Subcase &Group &Conditions\\\hline
	1 & 1a & $\Spin(6)$ & $n=3$\\\hline
      2 & 2a &$\Spin(2)\Spin(6)$ & $n=4$ \\\hline
      2 & 2b &$\Spin(6)\Spin(6)$ & $n=6$ \\\hline
      2 & 2c &$\Spin(2)\Spin(5)$ & $n=3$ \\\hline
      2 & 2d &$\Spin(5)\Spin(6)$ & $n=5$ \\\hline
      2 & 2e &$\hat\U$($n$) & -- \\\hline
      3 & 3a &$\hat\U(n-3)\Spin(6)$ & $n\geq 5$\\\hline
      3 & 3b &$\hat\U(n-2)\Spin(5)$ & $n\geq 4$\\\hline
      3 & 3c &$\Spin(2)\hat\U(n-1)$ & -- \\\hline
      4 & 4a &$\hat\U$($n_1$)$\hat\U$($n_2$) & $n_1, n_2\geq1, n_1+n_2=n$\\\hline
    \end{tabular}
\end{table}

\begin{itemize}
  \item[1a)] Transfering the diagram in to $\SO(7)$, we know that $(\H)_0$ is
    not $\SO(5)$ because the action of $\Spin(7)$ is effective. This implies
    $\H=\SU(3)$. We know $\K_-/\H$ is a sphere or a real projective space, but
    there is no embedding $\SU(4)\into\SO(7)$, so $\K_-=\G_2$. Since all
    embeddings $\G_2\into\SO(7)$ are conjugate, we can choose on such that
    $\SU(3)\subset\G_2$ acts trivially on $e_1$, which in turn completely
    determines $\K_+$, so we actually get one example here. It is primitive,
    because $\K_+$ is a maximal subgroup of maximal rank and there is no
    embedding $\K_-\into\K_+$.
  \item[2a)] We transfer the whole situation into $\SO(9)$ and have
    $\K_+=\S^1\SO(6)$ and $\H=\S^1\SU(3)$, where we can write
    $\S^1=\{\diag(1,z^l,z^k,\dots,z^k)\mid z\in\S^1\}$. But $\K_+/\H$ being a
    projective space implies $\SO(6)\cap\H=\SU(3)$, which can be written as
    
    \begin{displaymath}
      \SO(6)\cap\H = \{\diag(1,z^l,z^k,z^k,z^k)\cdot A\mid A\in\SU(3),z^l=1\}
    \end{displaymath}

    so that $z^{3k}=1$ whenever $z^l=1$, which implies that $l$ divides $3k$.
    But $\gcd(k,l)=1$, and by reparametrisation we can assume $l=1,3$. We now
    divide cases by $(\H_-)_0$:

    \begin{itemize}
      \item[$\bullet$] If $(\H_{-})_0=\SU(3)$, we know $\K_-/\H_-$ is $\U(2)$ or $\SO(3)$.
	But the normalizer $N(\SU(3))=\SO(3)\S^1$ does not contain $\SU(2)$, and
	$\K_-/\H_-=\SO(3)$ implies $k=0$, in which case we have that
	$\H\into\SO(9)$ is surjective on the fundamental group, a contradiction.
      \item[$\bullet$] If $(\H_{-})_0=\S^1$, we have $N(\S^1)=\S^1\U(3)$ for
	$k\ne 1$, so we will not get a primitive example in this case. If $k=1$,
	we have $N(\S^1)=\U(4)$, but both $\K_\pm$ are contained in $\SO(8)$, so
	again no primitive example can be found here.
      \item[$\bullet$] If $\H_-$ is finite, we have $\K_-=\S^1\SU(4)$ or
	$\K_-=\S^1\SO(6)$, but in both cases $\K_\pm\subset\SO(8)$ again, a
	contradiction to primitivity.
    \end{itemize}

  \item[2b)] Again, by virtue of $-1\in\K_+\setminus\H$ we look at the situation
    in $\SO(13)$. We have $\SO(6)\subset\H_+$, and therefore $l_\pm>1$, so that
    all of $\K_\pm,\H$ are connected. But $\SO(6)\subset\H$ implies that
    $\SO(13)/\H$ is simply connected, a contradiction, because
    $\Spin(13)/(\pi^{-1}(\H))_0$ is a nontrivial cover.

  \item[2c)] We transfer the discussion to $\SO(7)$. Then we have
    $\H=\S^1\SU(2)$, where $\S^1=\{\diag(z^l,z^k,z^k,1)\mid z\in\S^1\}$ (where
    $(k,l)=1$) and $\SU(2)$ is given accordingly. We have $\SO(5)/\SU(2)=\R P^7$
    so that $\SO(5)\cap\H=\SU(2)$, so that $z^{2k}=1$ whenever $z^l=1$, implying
    that $l$ divides $2k$. But $(k,l)=1$, so that $l$ divides 2, and therefore
    $l=\pm1,\pm2$. But the image of $\pi_1(\S^1)\to\pi_1(\SO(7))$ is given by
    $l\mod 2\in\Z_2$, implying $l=2$ (if $l=\pm1$, $\SO(7)/\H$ is simply connected,
    a contradiction as before) and $k$ is odd. We assume $l=2$ (possibly
    reparametrizing $\S^1$)and  divide cases by $(\H_-)_0$:

    \begin{itemize}
      \item[$\bullet$] If $\H_-=\H$, we have $\K_-\subset N(\H)$ and
	$N(\H)_0=\T^2\SU(2)\subset\K_+$, a contradiction to primitivity.
      \item[$\bullet$] If $(\H_-)_0=\S^1$, we have
	$N(\S^1)=\T^2\SU(2)\subset\K_+$ as before, a contradiction.  
      \item[$\bullet$] If $(\H_{-})_0=\SU(2)$, we have
	$N(\SU(2))=\S^1\SU(2)\SO(3)$. Since $\H/\H_-=\S^1$, we have
	$\K_-/\H_-=\U(2)$ or $\SO(3)$. But the only $\SU(2)$ in $N(\H_-)$ is
	also in $\H_-$, and the $\SO(3)$-factor does not contain $\S^1$, so
	there's no example in this case either.
      \item[$\bullet$] If $\H_-$ is finite, we have either $\K_-\simeq\U(3)$ or
	$\K_-\simeq\SO(2)\SO(5)$, where the second case of course does not give
	a primitive manifold. $\SU(3)$ is determined by $\H$, and the center of
	$\K_-$ is determined by $\H$ up to the first coordinate. It is either
	given by $\{\diag(z,z,z)\}$ or $\{\diag(\bar z,z,z\}$. But conjugation
	of the first coordinate fixes both $\K_+$ and $\H$, so we can assume
	$\K_-$ is the standard upper left block. Since $\U(3)\into\SO(7)$ is
	surjective on the fundamental group, $\pi^{-1}(\U(3))$ is connected and
	contains $-1$, so that $\K_-/\H$ is a projective space. Now $\SU(3)\cap
	\H=\SU(2)$ , which implies $k=1$ or $k=-3$ (note that
	$\SU(3)\cap\H/\SU(2)=\Z_2$ and that $z=\pm 1$ gives lements of
	$\SU(2)$). Those are obviously two primitive examples.
    \end{itemize}

  \item[2d)] This follows the exact same reasoning of case 2b above, and does
    not give a primitive example.
  \item[2e)] Since $-1\in\K_+$, again we transfer the situation to $\SO(2n+1)$.
    We have either $\SU(n)\not\subset\H_+$ or $n=3$ and $\H=\SU(3),
    \K_-\simeq\G_2$, or $n=4$ and $\H=\SU(4), \K_-\simeq\Spin(7)$, since
    otherwise $\SU(n+1)\subset\K_-$, but there is no embedding
    $\SU(n+1)\into\Spin(2n+1)$. The last case is not primitive, because
    obviously $\K_\pm\subset\SO(8)$. The second case gives one primitive example
    ($\K_+$ is maximal and does not contain a group isomorphic to $\K_-$, so it
    is primitive; since $\H$ uniquely determines $\K_+$, there is only one).
    Note that since $l_+=1$, components may occur, but the normalizer of $\G_2$
    is just the nontrivial element that maps to the identity in $\SO(7)$, so
    this will just produce an example of $\SO(7)$.
    
    From now on assume $\H=\S^1_k\SU(n-1)$ where

    \begin{displaymath}
      \S^1_k=\{(\diag(1,z^{(k+2)(n-1)},z^{-k},\dots,z^{-k})\mid z\in\S^1\}
    \end{displaymath}

    and $k$ is odd (see section \ref{subsubsec:u(n)onproj}). This implies
    $N(\H)_0=\T^2\SU(n-1)\subset\K_+$, so that $\H\ne\H_-$ and $l_{-}>1$. Also,
    $N(\SU(n-1))_0=\SO(3)\U(n-1)$, but the $\SO(3)$-factor does not contain
    $\S^1_k$, so in fact $l_->2$. Looking at both normalizers and checking
    against the list of effective transitive sphere actions
    (\ref{table:sphereactions}) we also see $l_->3$. For $n\geq 6$, this implies
    $\K_-\simeq\U(n)$, but it is already determined by $\S^1_k$ (note that $k$
    cannot be -2), so that $\K_\pm\subset\SO(2n)$ and no primitive example
    arises in this case.

    If $n\leq 5$, there are four additional possibilities for $\K_-$, two coming
    from the accidental isomorphisms $\Sp(2)\simeq \Spin(5)$ and
    $\SU(4)\simeq\Spin(6)$. They have already been considered in 2c and 2a
    respectively. For $n=4$, there is the possibility that $\K_-\simeq\S^1\G_2$.
    This is not possible since $k\ne 0$ implies that $\S^1$ is not in the
    normalizer of $\G_2$. For $n=5$ we can have $\K_-\simeq\S^1\Spin(7)$, but
    $k$ is odd so $\S^1$ does not commute with $\Spin(7)\subset\SO(8)$.

  \item[3a)] We have $-1\in\K_+$, so we transfer the situation to $\SO(2n+1)$
    again. We cannot have $\SO(6)\subset\H$, for that would imply that
    $\SO(2n+1)/\H$ is simply connected. But then $\SU(n-3)\subset\H_+$, because
    $\SO(6)\times\SU(n-3)$ can not act transitively almost effectively on a
    projective space. Now also $\SO(5)\not\subset\H$ for the same reasons as above, so we
    actually have $\H=\SU(3)\S^1\SU(n-3)$. Since $\K_-$ can not contain
    $\SU(3)\SU(n-2)$ (which has no embedding into $\SO(2n+1)$), we have that
    $\SU(n-3)$ is contained in the isotropy group of an almost effective
    tansitive action on a sphere of a group that does not contain $\SU(n-2)$ as
    a transitively acting factor. This implies $n=5$ or $n=6$. 

    If $n=5$, we have $\K_-=\SU(3)\S^1\SO(5)$ (note that $\SO(5)/\SU(2)$ is a
    projective space and $\Spin(5)=\Sp(2)$), but $\SO(5)$ is determined by
    $\SU(3)$ and $\SU(2)$, and it's easy to see $\K_\pm\subset\SO(6)\SO(5)$, so
    that this is not primitive.

    If $n=6$, we have $\H=\SU(3)\S^1\SU(3)$ and the above reasoning implies
    $\K_-=\SU(3)\S^1\SO(6)$, but again $\K_\pm\subset\SO(6)\SO(6)$.

  \item[3b)] Again, we look at the image of the diagram in $\SO(2n+1)$ and have
    $\K_+=\U(n-2)\SO(5)$. We have $\SO(5)\not\subset\H$ (because $\G/\H$ is not
    simply connected), so $\SU(n-2)\subset\H_-$, which implies
    $\H=\SU(n-2)\tilde \S^1\SU(2)$. Now $\SU(n-1)\SU(2)\subset\K_-$ is a
    contradiction (we can see there is no such embedding into $\SO(2n+1)$ by
    looking at the dimensions of representations), and since $\SU(n-2)\subset
    \H_-$ we have $n=4$ or $n=5$. But $n=5$ implies
    $\K_-\simeq\SU(2)\S^1\SO(6)$, which does not give a primitive example by
    3a).

    So we are left with $n=4$, and $\K_-=\SU(2)\S^1\SO(5)$. We know $\tilde
    \S^1\subset\H$ commutes with both $\SU(2)$-factors, so it is of the form

    \begin{displaymath}
      \tilde\S^1=\{\diag(z^{l_1},z^{l_1},1,z^{l_2},z^{l_2})\mid z\in\S^1\}
    \end{displaymath}

    for some $l_1,l_2\in\Z$ with $(l_1,l_2)=1$. Since complex conjugation on
    $\U(2)$ can be realized by a matrix in $\SO(4)$, we can reconjugate the
    diagram so that $l_1,l_2\geq 0$, without changing $\K_\pm$. Now since
    $\K_+/\H$ is a projective space, we know $\SO(5)\cap\H=\SU(2)$, so we know
    that $z^{l_1}=\pm1$ for all $z$ that satisfy $z^{l_2}=1$, i.e.\ $l_2\mid
    2l_1$. But $l_1,l_2$ are coprime, so $l_2\mid 2$, so $l_2=1$ or $l_2=2$. By
    symmetry reasons ($\K_-$ is a conjugate of $\K_+$ by just exchanging the
    first 4 coordinates with the last 4), also $l_1=1$ or $l_1=2$. Exchanging
    $\K_+$ and $\K_-$ does not change the manifold by section
    \ref{subsec:uniq}, we can assume $l_1\geq l_2$ and are left with 2
    examples. Since exchanging the first 4 coordinates with the last 4 is not in
    in the unity component of the normalizer of $\H$ (all automorphisms in
    $N(\H)_0$ are inner), those 2 examples are primitive.

  \item[3c)] Since $-1\in\K_+$, we translate the situation to $\SO(2n+1)$.
    First assume $\SU(n-1)\subset\H_+$. This implies $\K_-\simeq\U(n)$, so no
    primitive example arises. Note that $\H/\H_-$ contains $\S^1\SU(n-1)$, which
    also excludes the cases $n=3$ and $\K_-\supset\Sp(2)$ and $n=4$ and
    $\K_-\supset\SO(6)$.

    So now assume $\H=\T^2\SU(n-2)$. We have $\S^1\subset\H/\H_-$, so either
    $\H/\H_-=\S^1$ or $\H/\H_-=\S^1\SU(n-2)$. The latter case just implies
    $\K_-\simeq\K_+$ , which does not give a primitive example for the following
    reason: Note that if all of $\H,\K_\pm$ have a trivial one-dimensional
    subrepresentation from their inclusion into $\SO(2n+1)$, these are
    neccessarily the same, so $\K_\pm\subset\SO(2n)$. Otherwise $\H$ has three
    one-dimensional trivial subrepresentations, and then $\SO(3)\subset
    N(\H)_0$, which includes the change of coordinate which allows us to assume
    $\K_\pm$ act trivial on the same one-dimensional subspace, and
    $\K_\pm\subset\SO(2n)$ as before. 
    
    So we know $\H/\H_-=\S^1$, and we divide cases by $l_-=3$ and $l_{-}=2$:

    \begin{itemize}
      \item[$\bullet$] If $l_-=2$, then $\K_-/\H_-$ is either $\SO(3)$ or
	$\SU(2)$. The first case would imply that
	$\Spin(3)\subset(\pi^{-1}(\K_-))_0$ acts on the sphere $\S^{l_{-}}$, but
	$-1\in\Spin(3)$ has to acts trivial, a contradiction. In the second case
	$\K_-$ has a trivial 1-dimensional subrepresentation as a subgroup of
	$\SO(2n+1)$, which is necessarily given by $\R e_1$ (because we have
	$\SU(2)\subset N(\S^1\SU(n-2))$, so this is not a primitive example.
      \item[$\bullet$] If $l_-=3$, we have $\K_-/\H_-=\U(2)$, so that
	$\K_-=\T^2\SU(2)\SU(n-2)$. First assume that $\H$ as a subgroup of
	$\SO(2n+1)$ has exactly one 1-dimensional trivial subrepresentation.
	This is neccessarily $\R e_1$, and the same is true for $\K_-$, which
	also has exactly one 1-dimensional trivial subrepresentation on which
	$\H$ need to act trivially, too. The other possibility is that $\H$ has
	three 1-dimensional trivial subrepresentations, but then
	$N(\H)_0\supset\SO(3)\SU(n-2)$ and so by conjugation in $N(\H)_0$ we can
	assume that $\K_-$ acts trivially on $\R e_1$, again a contradiction to
	primitivity.

    \end{itemize}

  \item[4a)] Just as in 4a) of case $\G=\SO(2n+1)$, no primitive example can
    arise in this case. The proof carries over verbatim by transferring the
    situation to $\SO(2n+1)$.
\end{itemize}

\section{$\G=\Sp(n)$, $n\geq 2$}\label{sec:spn}

We claim that for $n>1$ the simply connected primitive cohomogeneity one
$\Sp(n)$-manifolds with positive Euler characteristic are given by table
\ref{table:spn}. Note that we will not treat $\Sp(2)\simeq\Spin(5)$ the way we
deal with the spin groups as described in section \ref{subsec:procspin}, but
just follow the procedure for the non-$\Spin$ groups.

\begin{table}[ht!]\caption{$\Sp(n)$-cohomogeneity one manifolds for
  $n>1$}\label{table:spn}
  \begin{center}
  \begin{tabular}{|c|} \hline 
    $\Sp(n_1-1)\Sp(n_2)\subset\Sp(n_1-1)\Sp(n_2+1),\Sp(n_1)\Sp(n_2)$\\\hline
    $\Sp(n-2)\Delta\Sp(1)\subset\Sp(n-2)\SO(2)\Delta\Sp(1),\Sp(n-1)\Sp(1)$\\
    where $\Delta\Sp(1)$ is the diagonal $\Sp(1)$ in the upper left
    $\Sp(2)$-block\\
    and $\SO(2)$ is the standard $\SO(2)$ in this same block\\\hline
    $\Z_2\Delta\Sp(1)\subset\SO(2)\Delta\Sp(1),\Z_2\Sp(1)\Sp(1), n=2$\\
    where $\Z_2$ is generated by $\begin{pmatrix}0 & 1\\ 1 &0\end{pmatrix}$
      and\\
    $\SO(2)$ and $\Delta\Sp(1)$ are as above\\\hline
    $\S^1\Sp(n-2)\subset\U(2)\Sp(n-2),\S^1\Sp(n-1)$\\
    where $\S^1=\{\diag(z,\bar z^2,1,\dots,1)\}$ or $\{\diag(z,1,\dots,1)\}$\\\hline
    $\U(n_1)\Sp(n_2-1)\subset\U(n_1+1)\Sp(n_2-1),\U(n_1)\Sp(n_2),
    n_1+n_2=n$\\\hline
    $\S^1_\H\Sp(n-2)\subset\Sp(1)\Sp(n-2),\S^1\Sp(n-1)$\\
    where $\S^1_\H=\{\diag(z,z^3,1,\dots,1)\}$ and $\Sp(1)\into\Sp(2)$ (upper left
    block) \\
    given by the irreducible $\SO(3)\into \SO(5)$\\\hline
    $\T^2_\H\Sp(n-3)\subset\S^1\SO'(3)\Sp(n-3),\T^2\Sp(n-2)$\\     
    where $\S^1=\{\diag(z,z,z,1,\dots,1)\}$, $\T^2_\H=\{\diag(1,z,\bar
    z,1,\dots,1)\}\cdot \S^1$\\
    and $\SO'(3)$ is the conjugation of the standard upper left \\$\SO(3)$-block
    by $\begin{pmatrix}1 & 0 & 0 \\ 0 & 1 & i\\0& 1&-i\end{pmatrix}$\\\hline
        \end{tabular}
  \end{center}
\end{table} 

The possibilities for $\K_+$ are given by table \ref{table:spnk+}.

\begin{table}[h]\caption{Possibilities for $K_+$}
  \label{table:spnk+}
  \begin{center}
    \begin{tabular}{|c|l|l|}\hline
      Factors& Group &Conditions\\\hline
      2 & $\U(n)$ & -- \\\hline
      2 & $\Sp(n_1)\Sp(n_2)$ & $n_1+n_2=n $\\\hline
      2 & $\S^1\Sp(n-1)$ & -- \\\hline
      2 & $\T^2$ & $n=2$ \\\hline
      3 & $\U(n_1)\Sp(n_2)$ & $n_1>1, n_1+n_2=n$\\\hline
      3 & $\S^1\U(n-1)$ & $n\geq3$\\\hline
      3 & $\S^1\Sp(n_1)\Sp(n_2)$ & $n_1 + n_2 + 1 =n$ \\\hline
      3 & $\Sp(1)\Sp(n_1)\Sp(n_2)$ & $n_1 +n_2+1=n$ \\\hline
      3 & $\T^2\Sp(n-2)$ & $n\geq 3$\\\hline
      4 & $\U(n_1)\U(n_2)$ & $n_1, n_2>1, n_1+n_2=n$\\\hline
      4 & $\Sp(1)\Sp(n_1)\Sp(1)\Sp(n_2)$ & $n_1,n_2\geq 1, n\geq 4$\\\hline
      4 & $\Sp(1)\Sp(n_1)\S^1\Sp(n_2)$ & $n_1,n_2\geq 1,n\geq 4$\\\hline
      4 & $\T^2\Sp(n_1)\Sp(n_2)$ & $n_1,n_2\geq 1,n\geq 4$\\\hline
      4 & $\Sp(1)\Sp(n_1)\U(n_2)$ & $n_1\geq 1,n_2\geq 2,n\geq 4$\\\hline
      4 & $\S^1\Sp(n_1)\U(n_2)$ & $n_1\geq 1,n_2\geq 2,n\geq 4$\\\hline
    \end{tabular}
  \end{center}
\end{table}

\begin{itemize} 
 \vskip 1cm 
  \item[2a)] If $n\ne 2$, $\SU(n)$ cannot be contained in the kernel, for that would mean
    $\SU(n)\subset\K_-\cap \H$ would act almost effectively on $\K_-/\H$, so
    that $\SU(n+1)\subset\K_-$. But there's no embedding of
    $\SU(n+1)\into\Sp(n)$ save the case $n=1$, which is not included here. This
    is easily seen by looking at the dimensions of representations.

    So assume $\SU(n)\not\subset\H_+$. Now we know $\H$ is given by
    $\S^1_k\SU(n-1)$ where $\S^1_k$ is given by

    \begin{displaymath} \S^1_k=\{\diag(z^{(k+1)(n-1)},z^{-k},\ldots,z^{-k})\mid
      z\in\S^1\} \end{displaymath}
    
    The unity component of the normalizer of $\H$ in $\Sp(n)$ is given by either
    $\Sp(1)\U(n-1)$ (for $k=-1$) or $\S^1\S^1_k\SU(n-1)$ (for $k\ne-1$), where
    $\S^1$ is the standard one in the left upper $\Sp(1)$ block. Let's first
    assume $\K_-\ne\U(n)$. For $k=-1$, we can have $\K_-=\Sp(1)\U(n-1)$, giving
    a primitve manifold. There's also the case $\K_-=\S^1\H$, where by
    conjugation in $N(\H)_0$ we can assume $\S^1\subset\K_-$ is standard. But
    this is not primitive, for $\K_-\subset\K_+$. For $k\ne -1$ there's only the
    case $\S^1\H$ where $\S^1$ is standard, which isn't primive.

    If $\K_-\simeq\U(n)$, we can determine the  embedding via the center as in case
    2a) of section \ref{sec:so2n}, which also leads to $\K_-=\bar\U(n)$ as the
    only primitive possibility, which only works in the cases $k=-1$ or $k=0$.
    But $\hat \U(n)=\sigma(\U(n))$ where $\sigma$ is conjugation by
    $\diag(j,1,\dots,1)\in N(\H)_0$ for $k=-1,0$, so no primitive manifold
    arises.

    Lastly assume $\SU(n)\subset\H_+$ which implies $n=2$ as said above. Since
    now $\H_0=\SU(2)$ and there's no embedding of $\SU(3)$ into $\Sp(2)$, we
    have $(\K_{-})_0\simeq\Sp(1)\Sp(1)$. We reconjugate the diagram so that
    $(\K_-)_0$ is standard, and $(\H)_0=\Delta \Sp(1)$ is diagonal. This already
    determines the maximal subgroup of maximal rank $\K_+=\S^1\Delta\Sp(1)$. If
    all groups are connected, this gives a primitive manifold. The normalizer of
    $(\K_{-})_0$ in $\Sp(2)$ is $\Z_2\Sp(1)\Sp(1)$, where the outer automorphism
    in $\Z_2$ switches the two $\Sp(1)$-factors. This gives $\H=\Z_2\Delta
    \Sp(1)$ and $\K_+$ as before, another primitive manifold.
  \item[2b)] First assume $\Sp(n_2)\subset\H_+$. If $n_1>2$, so that
    $\H=\Sp(n_1-1)\Sp(n_2)$ and we know from primitivity that
    $\Sp(n_2)\subset\H/\H_-$ and therefore $\K_-=\Sp(n_1-1)\Sp(n_2+1)$,
    resulting in a Grassmanian. Note that in the case $n_2=2$ we have
    $\H\simeq\Sp(n-3)\Spin(5)$, but we cannot have $\H\simeq\Sp(n-3)\Spin(6)$,
    because that would be a maximal rank subgroup, but not isomorphic to one of
    the groups we have given in table \ref{table:maxranksubgroups}. In the case
    $n_2=1$ we have $\H=\Sp(n-2)\Sp(1)$ so that we can have
    $\K_-=\Sp(n-2)\Spin(4)=\Sp(n-2)\Sp(1)\Sp(1)$ or $\K_-=\Sp(2)$ for $n=3$. For
    the first case, the standard representation of $\Sp(n)$ restricted to $\H$
    factors into 3 nonequivalent irreducible representations $\mathbb
    H^{n-2}\oplus\mathbb H\oplus\mathbb H$ where the last factor is acted on
    trivially. This already determines $\K_-$, which is standard, which would
    imply $\Sp(1)\subset\H_-\cap\H_+$, a contradiction. The latter case results
    in the action of $\Sp(3)$ on the Cayley plane.

    If $n_1\leq 2$ the possibilities given above arise as well, but we need to
    note that because of $\text{corank}(\H)=1$ by corollary \ref{cor:corH} we
    cannot have $\K_+/\H_+\simeq \SO(3)$ or $\SO(5)$.

    Of course, if $\Sp(n_1)\subset\H_+$ we may just exchange $\Sp(n_1)$ with
    $\Sp(n_2)$, so now we assume that none of both is contained in $\H_+$, which
    amounts to either $n_1=1$ or $n_2=1$ (just assume the latter) and
    $\H=\Sp(n-2)\Delta\Sp(1)$. The identity component of the normalizer of $\H$
    in $\G$ is given by $N(\H)_0=\S^1\Delta\Sp(1)\Sp(n-2)$, where $\S^1$ is
    given by $\SO(2)\into\Sp(2)$.  We can go through the possibilities for
    $\K_-$ by studying the isotropy representation. Two of the possibilities
    amount to $\K_-\simeq\K_+$, one of which is obviously not primitive
    ($\K_-=\K_+$), and the other one where $\K_-$ is obtained from $\K_+$ by
    exchanging the first two coordinates. This is given by conjugation with the
    matrix having

    \begin{displaymath} \begin{pmatrix}0&-1\\ 1&0\end{pmatrix} \end{displaymath}

    as the upper left 2x2-block, and extended by the identity to the other
    coordinates. This is in $N(\H)_0$ though, so no new primitive manifold
    arises. We can also have $\K_-=\Sp(1)\Sp(1)\Sp(n-2)$, where both possible
    cases (the additional $\Sp(1)$ can be in the diagonal or offdiagonal in the
    upper right block of $\Sp(2)$) are conjugate (in $\Sp(2)$ this conjugation
    is given by

    \begin{displaymath} 
      \frac{1}{\sqrt2}\begin{pmatrix} 1 & 1\\ -1 &
      1\end{pmatrix} 
    \end{displaymath}

  But this is in $N(\H)_0$ as well, so again no primitive example arises.  The
  last case is $\K_-=\S^1\Delta\Sp(1)\Sp(n-2)$ where $\S^1$ is given by the
  standard  embedding $\SO(2)\into\Sp(2)$, i.e.\ $\K_-\simeq\U(2)\Sp(n-2)$,
  which is primitive for the same reasons as above. Note that we have $l_-=1$
  here, so that components may occur. But the Normalizer of $\K_+$ in $\G$ is
  $\K_+$ itself save the case $n=2$, which was already treated in 2a.
    
\item[2c)] We cannot have $\K_+/\H_+=\S^1$, since that would imply
  $\Sp(n-1)\subset \H_+$ and therefore $\K_-=\Sp(n)$.

  The almost effective actions of $\S^1\Sp(n-1)$ on spheres
  are given by $(z,A)\cdot v = Avz^l$ with isotropy group
  $\Sp(n-2)\Delta\S^1_l$, where $\Delta\S^1_l$ is given by

  \begin{displaymath}
    \S^1_l:=\{\diag(z,z^l,1,\ldots,1)\mid z\in\S^1\}
  \end{displaymath}

  First assume $\Sp(n-2)\not\subset\H_-$. If $\Delta\S^1_l\subset\H_-$, we have
  $\K_-=\S^1_l\Sp(n-1)$ where $\S^1_l$ normalizes the $\Sp(n-1)$ factor. That is
  only possible for $l=0$ and $\K_-=\K_+$. If $\Delta\S^1_l\not\subset\H_-$,
  then we have $\K_-\simeq\K_+$ with only two choices for $\K_-$. One
  is $\K_-=\K_+$, and the other is $\K_-=\sigma(\K_+)$, where $\sigma$ is the
  coordinate change which exchanges the first 2 coordinates. But this implies
  $l=1$, where we have $N(\H)_0=\U(2)\Sp(n-2)$, which includes $\sigma$, so this
  manifold is not primitive as well.

  From now on we assume $\Sp(n-2)\subset\H_-$, that is, $\H/\H_-=\S^1$. First
  also assume $\K_-/\H_-=\S^2$, which implies $\K_-=\Sp(1)\Sp(n-2)$. Note that
  $\Sp(1)$ and $\Delta\S^1_l$ have to share a common central element of order 2.
  If $l$ is even, that element is $\diag(-1,1,\dots,1)$, which already implies
  that $\Sp(1)$ is standard in the upper left corner and $l=0$, obiously not a
  primitive manifold. If $n$ is odd, that element is $\diag(-1,-1,1,\dots,1)$,
  which is central in the upper left $\Sp(2)$ block, so we can look for
  embeddings $\Sp(1)\into\Sp(2)$ via embeddings $\SO(3)\into\SO(5)$. There are 2
  of those: The standard embedding, which corresponds to $l=1$ by the induced
  weights for the isotropy representation of $\S^1\simeq\SO(2)$. This will not
  be primitive, because $N(\H)_0$ contains $\U(2)$ in the upper left 2x2-block
  which acts transitively on the possible extensions $\S^1_1\into\Sp(1)$, and
  $\K_-=\Sp(n-2)\Delta\Sp(1)$ is contained in $\Sp(1)\Sp(n-1)$ as well as $\K_+$.
  The other possible embedding is given via the representation on the traceless
  symmetric 3x3-matrices by conjugation (in representation terms, that's given
  by (4)) and corresponds to $l=3$. Note that this representation is
  irreducible, so that $\K_-$ is not contained in $\Sp(1)\Sp(n-1)$, which is the
  maximal subgroup containing $\K_+$, so that the corresponding manifold is
  primitive.

  The last case is $\K_-/\H_-=\S^3$. Since, as before, $\Sp(n-2)\subset\H_-$, we
  shift the discussion into the upper left $\Sp(2)$-block. The isotropy
  representation of $\Delta\S^1_l$ in $\Sp(2)$ has 4 summands of dimension 2 and at
  least one trivial of dimension one. Only the two off-diagonal summands (of
  weight $l\pm1$) belong to groups not contained in $\Sp(1)\Sp(1)$, which would
  not lead to a primitive manifold, because $\K_-\subset\K_+$. The off-diagonal
  summands belong to the standard $\SU(2)$ in $\Sp(2)$ as well as its conjugate
  by

  \begin{displaymath}
    \begin{pmatrix}
      1& 0 \\ 0 & j
    \end{pmatrix}
  \end{displaymath}

  Note that this conjugation leaves $\K_+$ invariant and just exchanges
  $\Delta\S^1_l$ with $\Delta\S^1_{-l}$ in $\H$, so we may just assume
  $\K_-=\SU(2)\Delta\S^1_l$, which leaves $l=0,-2$ (see section
  \ref{subsubsec:un1spheres}). In both cases, $\K_\pm$ are two differen maximal
  subgroups of maximal rank, so the manifolds are primitive (for $l=0$ the
  manifold is homogeneous, a special case of 3a). 
\item[2d)] We have $\H_0=\S^1$, and by checking the isotropy representation as
  in the last paragraph of the previous case we can see $(\K_-)_0=\U(2)$ (note
  that again conjugation by $\diag(1,j)$ leaves $\K_+$ invariant, and $\H$
  additionally). Now $\K_+\subset \K_-$ and no primitive manifold arises.

\item[3a)] First, suppost $\Sp(n_2)\subset\H_+$, so that
  $\H=\S^1_k\SU(n_1-1)\Sp(n_2)$ where 
  
  \begin{displaymath}
    \S^1_k=\{\diag(z^{(k+1)(n_1-1)},\underbrace{z^{-k},\dots,z^{-k}}_{n_1-1},\underbrace{1,\dots,1}_{n_2})
    \mid z\in\S^1\} \end{displaymath} 

  By primitivity $\Sp(n_1)\not\subset\H_-$, so that we must have
  $\Sp(n_2+1)\subset \K_-$, which is contained in the normalizer of
  $\SU(n_1-1)$, which implies $k=-1$ if $n_1\ne 2$, which gives a
  primitive manifold. If $n_1=2$, we have $\K_-\simeq \S^1\Sp(n-1)$, which was
  treated in 2c.
  
  Now let $\SU(n_1)\subset\H_+$. Then $\K_+=\S^1\SU(n_1)\Sp(n_2)$ acts in the
  following way: For $(z,A)\in\S^1\Sp(n_2)$ we have $\S^{l_+}\ni v\mapsto
  Avz^{-kn_1}$. The isotropy group $\H$ is given by
  $\H=\S^1_k\SU(n_1)\Sp(n_2-1)$, where

  \begin{displaymath}
    \S^1_k=\{\diag(\underbrace{z,\dots,z}_{n_1},z^{kn_1},\underbrace{1,\dots,1}_{n_2-1})\mid
    z\in \S^1\}
  \end{displaymath}

  By primitivity $\SU(n_1)\not\subset\H_-$. For $n_1\ne 2$ this implies
  $\SU(n_1+1)\subset\K_-$, which implies $k=0$ and $\K_-=\U(n_1+1)\Sp(n_2+1)$ or
  its conjugate by

  \begin{displaymath}
    \diag(\underbrace{1,\dots,1}_{n_1},j,\underbrace{1,\dots,1}_{n_2})
  \end{displaymath}

  But this conjugation leaves $\H$ invariant as well as $\K_+$, so both
  manifolds are equivalent. Note that we cannot have $\S^1\subset\H_-$ for we do
  have $\S^1\subset\H_+$. This is one primitive example, the homogeneous space
  
  \begin{displaymath}
    \frac{\Sp(n+1)}{\U(n_1+1)\Sp(n_2)}
  \end{displaymath}

  If $n_1=2$, there is an additional posibbility:
  $\K_-\simeq\S^1\Sp(2)\Sp(n-3)$. We conjugate the whole diagram
  to make $\K_-$ standard. Then $\H=\S^1_l\Sp(1)\Sp(n-3)$ where

  \begin{displaymath}
    \S^1_l=\{\diag(z,z^l,1,\dots,1)\mid z\in\S^1\}
  \end{displaymath}

  and $\S^1_l\cap\Sp(1)$ is trivial. It is clear from this that
  $\K_+\simeq\S^1\Sp(1)\Sp(n-2)$, a contradiction to the original assumption
  $\K_+\simeq\U(2)\Sp(n-2)$ (note though that the case
  $\K_+\simeq\S^1\Sp(1)\Sp(n-2)$ will be treated later, and not give a primitive
  example).
  \item[3b)] We choose $\K_+$ such that $\SU(n-1)\subset\K_+$ is the lower right
  block. If we assume $\SU(n-1)\subset\H_+$, then by primitivity we have
  $\K_-\simeq\U(n)$, which is treated in case 2a. So we can assume
  $\SU(n-1)\not\subset\H_+$, which implies $\H=\T^2\SU(n-2)$. If
  $\SU(n-2)\not\subset\H_-$, then $\K_-\simeq\S^1\U(n-1)$. By primitivity,
  $\SU(n-1)\subset\K_-$ is the upper left block. But $\K_\pm\subset\U(n)$ where
  $\U(n)$ is either standard or has center $\diag(\bar z,z,\dots,z)$.

  So from now on we can assume $\SU(n-2)\subset\H_-$, i.e.\ $\H/\H_-\simeq\S^1$.
  This implies $l_-=2$ or $l_-=3$, and we divide cases by $\K_{-}$:

  \begin{itemize}
    \item[$\bullet$] If $l_-=2$ and $\K_-=\SU(2)\U(n-2)$, we reconjugate the
      diagram to make $\K_-$ standard. Then $\H=\S^1_1\U(n-2)$ where $\U(n-2)$
      is the lower right block and

      \begin{displaymath} \S^1_1=\{\diag(z,\bar z,1,\dots,1)\mid z\in\S^1\}
      \end{displaymath}

      This way we see there are 2 choices for $\SU(n-1)\subset\K_+$, both of
      which are equivalent. Now we claim that the factor $\SU(2)\subset\K_-$ is
      uniquely determined by $\S^1_1$ up to conjugation in $N(\S^1_1)_0$, which
      finishes this case, because then $\K_\pm\subset\U(n)$ up to conjugation in
      $N(\H)_0$.

      To prove the claim we restrict the discussion to the upper left $\Sp(2)$
      block. We have $\Sp(2)\simeq \Spin(5)$, but $\SU(2)\subset\Sp(2)$ contains
      the central element, so it corresponds to $\SO(3)\subset\SO(5)$, and
      $\S^1\subset\SU(2)\subset\Sp(2)$ to $\SO(2)\subset\SO(3)\subset\SO(5)$.
      But $\SO(3)\subset\SO(5)$ is uniquely determined up to conjugation in
      $N(\SO(2))$. So there is no primitive example in this case.

    \item[$\bullet$] If $l_-=2$ and $\K_-=\Sp(1)\SU(n-2)\S^1$, we claim
      $\K_\pm\subset\Sp(1)\Sp(n-1)$. As above, we restrict the discussion to the
      upper left $\Sp(2)$ block to show that $\Sp(1)$ is uniquely determined by
      $\S^1$ it contains. In $\SO(5)$, $\SO(2)$ has either three 1-dimensional
      trivial subrepresentations, and the $\SU(2)$ containing it is unique in
      $N(\SO(2))_0$, or only one, but then the $\SU(2)$ containing it is already
      determined. This shows the claim if $\Sp(1)$ does not contain the central
      element of $\Sp(2)$. If $\Sp(1)$ contains the central element of $\Sp(2)$,
      the discussion above shows that is is also uniquely determined by the
      $\S^1$ it contains. Again, no primitive example arises in this case.

    \item[$\bullet$] If $l_-=3$ then either $\K_-=\U(2)\U(n-2)$ or
      $\K_-=\S^1\Sp(1)\U(n-2)$ we can use the same arguments as above to show
      that the manifold is not primitive.
  \end{itemize}

\item[3c)] Suppose $\S^1\Sp(n_1)\subset\H_{+}$, so that
  $\H=\S^1\Sp(n_1)\Sp(n_2-1)$. We have $\S^1\Sp(n_1)\not\subset\H_{-}$ by
  primitivity, so that $\K_+\simeq \S^1\Sp(n_1+1)\Sp(n_2-1)$. But the $\S^1$
  factor of $\H$ is not diagonally embedded into the factors $\S^1\Sp(n_1+1)$ so
  that $\S^1\subset\K_-$ can only act trivially on $\S^{l_-}$, which implies
  $\S^1\subset\H_{-}$, a contradiction.

  $\Sp(n_1)\Sp(n_2)$ can only occur as the isotropy group of a transitive
  effective action on a sphere in the case $n_1=1$. But even in that case
  $\H=\Z_k\Sp(n_1)\Sp(n_2)$ cannot result in $\K_-=\Sp(1)\Sp(n_2+1)$ because
  as above $\Sp(1)\subset \H$ is not diagonally embedded into $\K_-$.

  So now we can assume $\H=\Sp(n_1-1)\Sp(n_2)\Delta\S^1_l$, where we assume
  $\Sp(n_1)$ in the upper left block, $\Sp(n_2)$ in the lower right, and

  \begin{displaymath}
    \Delta\S^1_l=\{\diag(\underbrace{1,\dots,1}_{n_1-1},z^l,z,
    \underbrace{1,\dots,1}_{n_2})\mid z\in\S^1\}
  \end{displaymath}

  Note that this choice actually means $\K_+=\Sp(n_1)\S^1\Sp(n_2)$ so that
  $\S^1$ acts on the $(n_1+1)$st coordinate.  If $|l|\ne1$, this implies
  $\K_-=\Sp(n_1-1)\S^1\sigma(\Sp(n_2+1))$, where $\sigma$ exchanges the
  Koordinates $n_1$ and $n_1+1$. This is not primitive, for
  $\K_\pm\subset\varphi(\S^1\Sp(n-1))$ where $\varphi$ exchanges the first with
  the $(n_1+1)$st coordinate. If $|l|=1$, we can actually have
  $\Sp(n_2+1)\subset\K_-$ embedded as the lower right block, but in this case
  it's even easier to see that the resulting manifold is not primitive.

\item[3d)] This case is completely analogous to 3c (for $l=1$) in showing that the only
  possibility is $\H=\Sp(n_1-1)\Sp(n_2)$ which does not lead to a primitive
  manifold again (because $\K_\pm\subset\Sp(1)\Sp(n-1)$).

\item[3e)] First suppose $\Sp(n-2)\subset\H_+$. This implies
  $\Sp(n-1)\subset\K_-$ and $\K_-$ is a subgroup of maximal rank, so we can
  refer to one of the previous cases.

  So we can assume $\H=\T^2\Sp(n-3)$. If $\Sp(n-3)\not\subset\H_-$, we have
  $\K_-\simeq\K_+$, and we claim that this is not a primitive manifold. For this
  note that both $\K_\pm$ (acting on $\mathbb H^n$) have 2 one-dimensional
  subrepresentations, while $\H$ has 3, which contain the ones of $\K_\pm$. So
  we see that $\K_\pm$ have a common one-dimensional subrepresentation, which
  shows $\K_\pm\subset\Sp(1)\Sp(n-1)$, so the manifold is not primitive.

  We are left with the case $\Sp(n-3)\subset\H_-$, which implies $l_{-}=2$ or
  $l_-=3$. We divide cases by $\K_+$

  \begin{itemize}
    \item[$\bullet$] If $l_-=2$ and $\K_-\simeq\SO(3)\S^1\Sp(n-3)$, we have
      $\S^1\subset N(\SO(3)\Sp(n-3)$. We reconjugate the diagram so that $\K_-$
      is standard and $\SO(2)\subset\H/\H_-$ is the upper left block. Then
      furthermore conjugate the diagram by the matrix having

      \begin{displaymath}
	\begin{pmatrix}1 & i\\ 1 & -i\end{pmatrix}
      \end{displaymath}

      in the upper left block, extended by the identity matrix $E_{n-2}$. This
      conjugates $\SO(2)\subset\H$ into 
      
      \begin{displaymath}
	\{\diag(z,\bar z,1,\dots,1)\mid z\in\S^1\}
      \end{displaymath}

      Now $\S^1\subset\H_-$ is in $N(\SO(3)\Sp(n-3))$ and therefore
      $\T^2\subset\H$ is given by

      \begin{displaymath}
	\{\diag(z,z,z,1,\dots,1)\mid z\in\S^1\}\cdot\{\diag(z,\bar
	z,1,\dots,1\}\mid z\in\S^1\}
      \end{displaymath}

      This completely determines $\K_+$ and we find one primitive example.
    \item[$\bullet$] If $l_-=2$ and $\K_-\simeq\SU(2)\S^1\Sp(n-3)$, we
      reconjugate the diagram to put $\K_-$ into standard form, which gives
      $\H=\S^1_1\S^1_2\Sp(n-3)$ where $\Sp(n-3)$ is the lower right block and

      \begin{align*}
	\S^1_1&=\{\diag(z,\bar z,1,\dots,1)\mid z\in\S^1\}\\
	\S^1_2&=\{\diag(1,1,z,1,\dots,1)\mid z\in\S^1\}
      \end{align*}

      From this it is clear that $\Sp(n-2)\subset\K_+$ is the lower right block.
      But $\H\cap\Sp(n-2)\ne\Sp(n-3)$, so $\K_+/\H$ is not a sphere.
    \item[$\bullet$] If $l_-=3$ and $\K_-\simeq\U(2)\S^1\Sp(n-3)$, we again
      reconjugate the diagram to make $\K_-$ standard. Then $N(\H)_0$ contains
      $\SO(3)$, so we can assume that $\Sp(n-2)\subset\K_+$ is the lower right
      block, so we have $\K_\pm\subset\Sp(2)\Sp(n-2)$, and the manifold is not
      primitive.
  \end{itemize}
\item[4a)] The fact that $\U(n)$ is given as the centralizer of its center
  $\S^1$ in $\Sp(n)$ leads to the same reasoning as in case 4a of
  \ref{sec:so2n} to show that no primitive manifolds arise in this case.

\item[4b-4d)] We will show that no primitive manifold can arise in these cases,
  because both $\K_\pm$ have a common one-dimensional subrepresentation when
  acting on $\mathbb H^n$, showing $\K_\pm\subset\Sp(1)\Sp(n-1)$. For this, denote the
  number of one-dimensional subrepresentations of the standard representation of
  a subgroup $\H\into\Sp(n)$ by $s(\H)$. Note that if $\K_\pm$ acts on a
  one-dimensional subspace of $\mathbb H^n$, so does $\H$.

  \begin{itemize}
    \item[$\bullet$] If $n_1,n_2\geq 2$, then $s(\K_+)=2,s(\H)=3,s(\K_-)=2$ and
      it is clear that both $\K_\pm$ share a one-dimensional invariant subspace.
    \item[$\bullet$] If $n_1=1$ and $\Sp(n_1)\not\subset\H_+$, then $s(\K_+)=3,
      s(\H)=3,s(\K_-)=2$ and it is clear that the claim holds.  
    \item[$\bullet$]
      If $n_1=1$, $\Sp(n_1)\subset\H_+$ and $n_2>1$, then
      $s(\K_+)=3,s(\H)=4,s(\K_-)=3$ and the claim holds.
    \item[$\bullet$] If $n_1=1=n_2$, then $s(\K_+)=4=s(\H),s(\K_-)=3$ and the
      claim holds.
  \end{itemize}

\item[4e)] Suppose $\H=\Delta\Sp(1)\Sp(n_1-1)\U(n_2)$. Then
  $\SU(n_2+1)\subset\K_-\cap N(\Delta\Sp(1)\Sp(n_1-1)$, a contradiction. So we
  can assume $\H=\Sp(1)\Sp(n_1)\U(n_2-1)_k$ for some $k\in\Z$. Then
  $\K_-=\Sp(1)\Sp(n_1+1)\U(n_2-1)$ and $k=1$. If $n_2\geq 2$, this is exactly
  the previous case by exchanging $\K_-$ and $\K_+$. If $n_2=2$, then
  $N(\H)\supset\SO(3)$, so we can assume $\Sp(n_1+1)\subset\K_{-}$ is the upper
  left block, which shows $\K_\pm\subset\Sp(n_1+1)\Sp(n_2)$.

\item[4f)] First assume $\H=\T^2\Sp(n_1-1)\SU(n_2)$. This implies
  $\K_-\simeq\T^2\Sp(n_1-1)\SU(n_2+1)$ by the classification of actions on a
  sphere (note that this is true even in the case $n_2=2$, because
  $\Sp(2)/\SU(2)$ is not a sphere). By primitivity, $\SU(n_2+1)\subset\K_-$ is
  completely determined for $n_1>1$ and we have $\K_\pm\subset\Sp(n_1)\Sp(n_2+1)$, so the
  manifold is not primitive. If $n_1=1$, we have $\SO(2)\subset N(\H)_0$, so we
  can assume that $\SU(n_2+1)\subset\K_-$ is the lower right block and the claim
  holds, too.

  The second possibility is $\H=\T^2\Sp(n_1)\SU(n_2-1)$. If
  $\K_-\simeq\T^2\Sp(n_1+1)\SU(n_2+1)$ and $n_2>2$ we can argue as before to show
  $\K_\pm\subset\Sp(n_1+1)\Sp(n_2)$. If in addition $n_2=2$, we have
  $\SO(3)\subset N(\H)_0$ so we can assume $\Sp(n-2)\subset\K_-$ is the upper
  left block, and the claim holds.

\end{itemize}
\section{$\G=\SO(2n)$, $n\geq 4$}\label{sec:so2n}

We claim that for $n\geq 4$ the simply connected primitive cohomogeneity one
$\SO(2n)$-manifolds with positive Euler characteristic are given by table
\ref{table:so2n}

\begin{table}[ht]\caption{$\SO(2n)$-cohomogeneity one manifolds for $n\geq
  4$}\label{table:so2n}
  \begin{center}
  \begin{tabular}{|c|} \hline 
    $\SU(4)\subset\Spin(7),\U(4), n=4$\\\hline
    $\S^1\SU(4)\subset\S^1\Spin(7),\U(5), n=5$\\\hline
    $\S^1\SU(n-1)\subset\sigma(\U(n)),\U(n)$\\
    where $\S^1=\{\diag(z,1,\dots,1)\}$ or $\S^1=\{\diag(1,z,\dots,z)\}$ and\\
    $\sigma$ is conjugation by $\diag(-1,1,\dots,1)$\\\hline
    $\SO(2n_1)\SO(2n_2-1)\subset\SO(2n_1+1)\SO(2n_2-1),\SO(2n_1)\SO(2n_2)$\\
    where $n_1+n_2=n$\\\hline
    $\Z_2\SO(2n-2)\subset\Z_2\SO(2n-1),\SO(2)\SO(2n-2)$\\
    where $\Z_2\subset\SO(2)$\\\hline
    $\U(n_1-1)\SO(2n_2)\subset\U(n_1-1)\SO(2n_2+1),\U(n_1)\SO(2n_2)$\\
    where $n_1+n_2=n$\\\hline
    $\T^2\SU(n-2)\subset\SO(3)\S^1\SU(n-2),\SO(2)\U(n-1)$\\
    where $\S^1=\{\diag(1,1,1,1,z,\dots,z)\}$, 
    $\T^2=\{\diag(z,1,\dots,1)\}\cdot \S^1$\\
    and $\SO(3)$ is the upper left block\\\hline
  \end{tabular}
  \end{center}
\end{table}

The possibilities for $\K_+$ are given by table \ref{table:so2nk+}.

\begin{table}[h]\caption{Possibilities for $\K_+$}
  \label{table:so2nk+}
  \begin{center}
    \begin{tabular}{|c|l|l|}\hline
      Factors& Group &Conditions\\\hline
      2 & $\U(n)$ & -- \\\hline
      2 & $\SO(2n_1)\SO(2n_2)$ & $n_1,n_2>2, n_1+n_2=n $\\\hline
      2 & $\SO(2)\SO(2n-2)$ & -- \\\hline
      3 & $\U(n_1)\SO(2n_2)$ & $n_1>1, n_2>2, n_1+n_2=n$\\\hline
      3 & $\U(n-1)\SO(2)$ & -- \\\hline
      3 & $\SO(2n-4)\SO(4)$ & $n>4$\\\hline
      4 & $\U(n_1)\U(n_2)$ & $n_1, n_2>1, n_1+n_2=n$\\\hline
      4 & $\U(n-2)\SO(4)$ & --\\\hline
      4 & $\SO(4)\SO(4)$ & $n=4$\\\hline
    \end{tabular}
  \end{center}
\end{table}

If a group of complex matrices is involved (e.g. $\U(n)$), we will deliberately
use complex notation for the corresponding real matrices. In particular, for
$e^{i\varphi}=z\in S^1$, we will use $\diag(z,\ldots,z)$ for the matrix
containing $\begin{pmatrix}\cos\varphi &\sin\varphi \\ -\sin\varphi
&\cos\varphi\end{pmatrix}$ on the diagonal $2\times2$-blocks and 0 everywhere
else.

\begin{itemize}
\item[2a)] If $\SU(n)\subset\H_+$, we have $n=4$ and $\K_-\simeq\Spin(7)$,
  because otherwise $\SU(n+1)\subset\K_-$, a contradiction to $\K_-$ being a
  subgroup of $\SO(2n)$. Assume the first case, then this gives one primitive
  example (it is primitive, because $\K_+$ is maximal and does not contain
  $\K_-$ or a conjugate; it is only one because $\K_+$ is uniquely determined by
  $\H$).
  
  So now we assume $\H$ is given by $\S^1_k\SU(n-1)$, where

  \begin{displaymath}
  \S^1_k=\{\diag(z^{(k+1)(n-1)},z^{-k},\dots,z^{-k})\}
  \end{displaymath} 

  by subsection \ref{subsubsec:un1spheres} of the appendix. We have
  $N(\SU(n-1))_0=\S^1\U(n-1)\subset\K_+$, therefore $\SU(n-1)\not\subset\H_-$
  and therefore $\K_-$ is isomorphic to $\U(n)$ save the case $n=5$ where
  $\K_-\simeq\S^1\Spin(7)$ is also possible. In the latter case it is clear that
  $k=0$ because $\S^1\subset\K_-$ must commute with $\Spin(7)\subset\SO(8)$.
  This is one primitive example, because $\Spin(7)$ is determined up to
  conjugation in $\SO(8)$ which also fixes $\H$. For primitivity note that
  $\K_-$ cannot be embedded into $\K_+$ and $\K_+$ is a maximal subgroup of
  maximal rank.

 From now on assume $\K_-\simeq\U(n)$. The center of $\K_-$ can be conjugated
 into the standard diagonal $\S^1$, and must of course commute with $\SU(n-1)$,
 so it is determined up to the first complex coordinate, which can be $z$ or
 $z^{-1}$. By primitivity, $\K_-\ne \K_+$, so the center of $\K_-$ is given by
 $z\mapsto\diag(z^{-1},z,\ldots,z)$, and therefore $\K_-=\varphi(\U(n))$, where
 $\varphi$ is conjugation by $\diag(-1,1,1,\ldots,1)$. Since $\varphi$ only
 fixes $\H$ for $k=0,-1$, we are left with 2 comomogeneity-1-manifolds, both of
 which are primitive indeed: Conjugation by
 
 \begin{displaymath}
   \diag(-1,1,\dots,1)
 \end{displaymath} 
 
 on $\U(n)$ has the same image as conjugation by

  \begin{displaymath}
    \diag(1,1,-1,1,\dots,-1,1)
  \end{displaymath} 
  
  (they only differ by $\diag(-1,1,\dots,-1,1)$, which is complex conjugation),
  but this restricts to complex conjugation on $\U(n-1)$, which is an outer
  automorphism and cannot possibly come from an element of
  $N(\H)_0=\S^1\U(n_1)$. Note that both of $\K_\pm$ must be connected, for
  $l_{\pm}=2n-1 > 1$.

\item[2b)] If $\K_+=\SO(2n_1)\SO(2n_2)$ with $n_1,n_2>2$, we may assume that $\H
  = \SO(2n_1)\SO(2n_2-1)$ is standard, and since $l_+=2n_1-1>1$ we know $\K_-$
  is connected. Since the action of $\SO(2n)$ is primitive and almost effective,
  we have  $\K_-=\SO(2n_1+1)\SO(2n_2-1)$. This is the primitive action of
  $\SO(2n)$ on the Grassmanian
  $\G_{2n_2}(\R^{2n+1})=\SO(2n+1)/\SO(2n_1+1)\SO(2n_2)$

\item[2c)] $\K_+/\H_+$ is not as $\O(2n-2)$ since $\O(2n-2)$ is not the direct
  product of $\SO(2n-2)$ and $\Z_2$. If $\K_+/\H_-$ is $\SO(2n-2)$, this is the
  exact same as case 2b), giving $\K_-=\SO(3)\SO(2n-3)$. So $\K_+/\H_-= \S^1$,
  which implies $\H=\Z_k\SO(2n-2)$ for some $k$. Now $\H_-$ can't contain
  $\SO(2n-2)$, so $\K_-$ must contain $\SO(2n-1)$ as a transitively acting
  factor.  Since the normalizer of $\SO(2n-1)$ in $\SO(2n)$ is
  $\S(\O(1)\O(2n-1))$ we are restricted to $k=1$ or $k=2$ (note that $l_+=1$, so
  actually components may occur). For $k=1$ this is analogous to case 2b), the
  action of $\SO(2n)$ on the Grassmanian $\SO(2n+1)/\SO(2n-1)\SO(2)$. For $k=2$
  we know $\SO(2n-1)\Z_2$ is isomorphic to $\O(2n-1)$ (since $\Z_2$ acts by an
  inner automorphism), but $\K_-/\H_-$ is not  $\O(2n-1)$ since this would imply
  $\H/\H_-\simeq \O(2n-2)$ which is not isomorphic to $\SO(2n-2)\times \Z_2$,
  which in turn $\H$ is.  So $\K_-/\H_-=\SO(2n-1)$, giving a non-effective
  action on $\mathbb CP^n$.
\item[3a)] $\K_+/\H_+$ is not $\SO(2n_2)$, since that would imply that $\K_-$
  has $\SU(n_1+1)$ as a simple factor as well as $\SO(2n_2-1)$. So $\K_+/\H_+
  =\U(n_1)$, and therefore $\H=\U(n_1-1)_k\SO(2n_2)$. Since the action is
  primitive and almost effective, $\K_-/\H_-$ must contain $\SO(2n_2+1)$ as a
  simple factor, and since $l_+=2n_1+1 > 1$,  $\K_-$ must be connected. So
  $\K_-=\U(n_1-1)_k\SO(2n_2+1)$, but the normalizer of $\U(n_1-1)_k$ does not
  contain $\SO(2n_2+1)$ save the case $k=0$. The embedding of $\K_-$ and $\H$
  into $\SO(2n)$ gives a faithfull representation, and we can see that there is
  a 2-dimensional trivial subrepresentation for $\H$, while there is a
  1-dimensional trivial subrepresentation for $\K_-$. It is clear that the
  latter is contained in the former and determines $\K_-$. But
  $N(\H)_0=\U(n_1-1)\SO(2)\SO(2n_2)$ acts transitively on the 2--dimensional
  trivial representation (via the factor $\SO(2)$), so by conjugating with an
  element of $N(\H)_0$ we can achieve $\K_-$ is standard.

This is the action of $\SO(2n)$ on
$\SO(2n+1)/\U(n_1)\SO(2n_2+1)$, which is primitive because $\K_+$ is a maximal
connected subgroup of $\SO(2n_1)\SO(2n_2)$ which in turn is maximal connected in
$\SO(2n)$. Both of these do not contain $\K_-$, even up to
components.

\item[3b)] In addition to the possibility from 3a) (giving the action of
  $\SO(2n)$ on $\SO(2n+1)/\U(n-1)\SO(3)$), there are 2 other possibilities.
  First, if $\H=\U(n-2)\SO(2)$, we can have $\K_-=\U(n-2)\U(2)$. From the
  isotropy representation of $\H$ in $\G$ we see that $\U(2)$ is contained in
  the lower right $\SO(4)$ block. It is determined by $\SU(2)\subset \U(2)$, and
  it is known that $\SO(4)=\SU(2)\SU(2)$, so there are 2 choices for the
  $\SU(2)$ factor. One of those leads to $\K_-$ being contained in $\U(n)$, not
  giving a primitive manifold (for $\K_+\subset\U(n)$ as well). The other one is
  obtained from the first by an outer automorphism, given by conjugation with
  $\diag(1,\ldots,1,-1)$. But this automorphism fixes both $\H$ and $\K_+$, so
  the resulting manifold is equivalent to the first one.  
  
  Secondly we can have $\K_+/\H\simeq\U(n-1)$, so that $\H=\T^2\SU(n-2)$. First
  assume $\SU(n-2)\not\subset\H_-$, which implies $\K_+\simeq\K_-$. Both
  $\K_\pm$ have a 2-dimensional subrepresentation from their embedding into
  $\SO(2n)$, which is neccessarily given by $\R e_1\oplus \R e_2$ or $\R
  e_3\oplus \R e_4$, and we chose $\K_+$ the way that this is $\R e_1\oplus\R e_2$. If
  it's the same for $\K_-$, we have $\K_\pm\subset\U(n)$, so the manifold is not
  primitive. If it is $\R e_3\oplus\R e_4$, we look at the center of
  $\U(n-1)\subset\K_-$. If it is given by

  \begin{displaymath}
    \{\diag(z,1,1,z,\dots,z)\mid z\in\S^1\}
  \end{displaymath}

  we again have $\K_\pm\subset\U(n)$. If it is given by

  \begin{displaymath}
    \{\diag(\bar z,1,1,z,\dots,z)\mid z\in\S^1\}
  \end{displaymath}

  we have $\K_\pm\subset\sigma(\U(n))$, where $\sigma$ is conjugation by
  $\diag(-1,1,\dots,1)$.

  So now we can assume $\SU(n-2)\subset\H_-$, which implies $\H/\H_-=\S^1$, and
  therefore $l_-=2$ or $l_-=3$. We divide cases by $\K_-$:

  \begin{itemize}
    \item[$\bullet$] If $\K_-=\S^1\SU(2)\SU(n-2)$, we look at the centralizer
      of $\SU(2)$ in the upper left $\SO(4)$-block: It is either given by
      $\diag(z,z)$, so that $\K_\pm\subset\U(n)$, or $\diag(\bar z,z)$, which
      implies $\K_\pm\subset\sigma(\U(n))$, where $\sigma$ is as above.
    \item[$\bullet$] If $\K_-=\SO(3)\S^1\SU(n-2)$ with $\SO(3)\subset
      N(\S^1\SU(n-2))$, we can reconjugate the diagram to achieve that $\K_-$
      acts trivially on $\R e_4$. This implies $\H=\S^1_1\S^1_2\SU(n-2)$ with
      $\SU(n-2)$ in the lower right block and

      \begin{align*}
	\S^1_1&=\{\diag(z,1,\dots,1)\mid z\in\S^1\}\\
	\S^1_2&=\{\diag(1,1,1,1,z,\dots,z)\mid z\in\S^1\}
      \end{align*}

      We note that $\SU(n-1)\subset\K_+$ is in $N(\S^1_1)$, so it is in the
      lower right block. This gives one primitive example.
    \item[$\bullet$] If $\K_-\simeq\U(2)\U(n-2)$, we can use the same reasoning
      above to see that the manifold is not primitive.
  \end{itemize}

\item[3c)]There are 2 actions arising in the same way as described in 2b),
namely $\SO(2n)$ acting on either $\SO(2n+1)/\SO(2n-3)\SO(4)$ or
$\SO(2n+1)/\SO(2n-4)\SO(5)$. We only need to note that $l_+$ is odd by corollary
\ref{cor:corH}.
\item[4a)] None of the actions arising are primitive. We have
\begin{displaymath}
\H=\T^2\SU(n_1)\SU(n_2-1)\text{ and }\H_+=\S^1\SU(n_1)
\end{displaymath}
Since $\H_-\cap\H_+$ is finite, $\SU(n_1)\subset\H$ cannot act trivial in
$\K_-$, and therefore $K_-$ is isomorphic to $\U(n_1+1)\U(n_2-1)$. As in case
2a) we can determine the embedding via the center. The standard embedding yields
$\K_-,\K_+\subset\U(n)$ and therefore not a primitive action, so the center has
to be given by
\begin{displaymath}
(z_1,z_2)\mapsto(\underbrace{z_1,\ldots,z_1}_{n_1},z_1^{-1},z_2,\ldots,z_2)
\end{displaymath}
which implies that bot $K_-$ and $K_+$ are contained in the subgroup of
$\SO(2n)$ isomorphic to $\U(n)$ with center
\begin{displaymath}
\left\{(\underbrace{z,\ldots,z}_{n_1},z^{-1},z^{-1},\ldots,z^{-1})\mid
z\in\S^1\right\}
\end{displaymath}

\item[4b)] This case is completely analogous to case 3a), because $\SO(4)$
cannot act as $\SO(3)$ for the same reason as in 3c).

\item[4c)] Again, no new example other than the ones arising from case 2b can
  occur because $l_+$ is odd by corollary \ref{cor:corH}.
\end{itemize}

\subsection{$\Spin(2n)$, $n\geq 4$}\label{subsec:spin2n}

We claim that there are no simply connected primitive cohomogeneity one
$\Spin(2n)$-manifolds with positive Euler characteristic for $n\geq 4$ such that
the element -1 (which projects to the identity, but is not the identity itself)
does not act trivially (i.e.\ it is not an action of $\SO(2n)$).

First, it is easily seen that in the case $n=4$, there can only be examples that
are coming from an $\SO(8)$-example. Consider the maximal torus $\T^4$ of
$\Spin(8)$, containing the maximal torus $\T^3$ of $\H$ (the corank of $\H$ is 1
by corollary \ref{cor:corH}). The group of involutions $\Z^4_2$ in $\T^4$
contains the center $\Z^2_2$ of $\Spin(8)$, which therefore has nonempty
intersection with the group of involutions $\Z^3_2$ of $\T^3$. Therefore, $\H$
contains an element of the center of $\Spin(8)$. Since the quotient of
$\Spin(8)$ by any central $\Z_2$ is $\SO(8)$, this shows the claim.

Now we have reduced the table of possibilities for $\K_+$ to those given in table \ref{table:spin2nk+}(where $\hat \U(n)$ is the preimage of $\U(n)\subset\SO(2n)$), where $n>4$.

\renewcommand\arraystretch{1.3}
\begin{table}[h]\caption{Possibilities for $K_+$}
  \label{table:spin2nk+}
  \centering
    \begin{tabular}{|c|l|l|}\hline
      Factors& Group &Conditions\\\hline
      2 & $\hat\U$($n$) & -- \\\hline
      3 & $\hat\U(n-1)\Spin(2)$ & -- \\\hline
      3 & $\hat\U(n-3)\Spin(6)$ & --\\\hline
      4 & $\hat\U$($n_1$)$\hat\U$($n_2$) & $n_1, n_2>1, n_1+n_2=n$\\\hline
    \end{tabular}
\end{table}

The last case is easily dismissed as well, since we already know from 
section \ref{sec:so2n} that no such manifold will be primitive (the proof
carries over almost verbatim).

\begin{itemize}
  
  \item[2a)]We look at the projection of the diagram in $\SO(2n)$. For the same
    reasons as in 2a) of section \ref{sec:so2n} we know $\SU(n)\not\subset \H$,
    so $\H=\S^1_k\SU(n-1)$, where

    \begin{displaymath}
     \S^1_k=\{\diag(z^{(k+2)(n-1)},z^{-k},\ldots ,z^{-k})\mid z\in\S^1\}
    \end{displaymath}
    
    and $k$ is odd. As already shown in 2a) of section \ref{sec:so2n} we have
    $\K_-\simeq \U(n)$ as well (note that the case $n=5$ and
    $\K_-\simeq\S^1\Spin(7)$ cannot occur here for it implies $k=0$), and since
    the diagram is primitive, $\K_-=\varphi(\U(n))$ where $\varphi$ is the
    conjugation by $\diag(-1,1,\ldots,1)$, which restricts to complex
    conjugation of the first component on $\S^1_k$. The latter is only invariant
    under $\varphi$ in the cases $k=0,-2$, but the quotients are spheres in
    these cases and so there are no new cases here.

  \item[3a)] As before, we will consider the situation in $\SO(2n)$. As in case
    3a) of section \ref{sec:so2n}, we have that $\H$ is given as
    $\U(n-2)_k\SO(2)$, where $\U(n-2)_k$ is determined by

    \begin{displaymath}\S^1_k:=\{\diag(1,1,z^{(k+2)(n-2)},z^{-k},\ldots,z^{-k})\mid
      z\in\S^1\}\end{displaymath} 
     
      with $k$ odd. This implies $\K_-/\H_-$ is given by either $\SO(3)$ or
      $\U(2)$. But both cases imply $k=-2$, a contradiction.    
    
   \item[3b)] Again, we will consider the diagram of subgroups of $\SO(2n)$. If
     $\K_+/\H_+=\SU(4)$, we have $\H=\SU(n-3)\S^1\SU(3)$, a contradiction,
     because $\K_-$ would have to contain $\SU(n-2)$, but there is no embedding
     of $\SU(n-2)\SU(3)$ into $\SO(2n)$. The same reasoning shows
     $\H\ne\U(n-3)\text{O}(5)$. We are left with $\H=\S^1_k\SU(n-4)\SO(6)$ where
     $\S^1_k$ is given by 

     \begin{displaymath}
       \S^1_k=\{\diag(z^{(k+2)(n-4)},z^{-k},\ldots,z^{-k}\underbrace{1,\ldots,1}_{\text{6
       times}})\}
     \end{displaymath}

     with $k$ odd. But in that case $\SO(2n)/\H$ is simply connected, a
     contradiction.\end{itemize}

\appendix 
\section{Appendix}
\subsection{The result}\label{subsec:result}
The entries of the following tables are of the form $\H\subset\K_-,\K_+$, where
$\K_+$ has maximal rank in $\G$, and $\K_\pm/\H$ are spheres. The manifold $M$
is then equivalent to

\begin{displaymath}
  \G\times_{\K_-}D^{l_-+1}\cup\G\times_{\K_+}D^{l_++1} 
\end{displaymath}

where $D^{l_\pm+1}$ is the unit ball with boundary $\partial
D^{l_-\pm+1}=\K_{\pm}/\H$, where the action of $\K_\pm$ is extended linearly
from the sphere $\K_\pm/\H$, and $\G\times_{\K_\pm} D^{l_\pm+1}$ is the quotient of
$\G\times D^{l_\pm+1}$ by the diagonal action of $\K_\pm$ (see
\ref{subsec:cohom1actions} for more details).

The tables for the $\Spin$-groups only list the effective actions (the
non-effective ones are listed for the $\SO$-groups), but those are given as
$\pi(\H)\subset\pi(\K_-),\pi(\K_+)$, where $\pi:\Spin(n)\to\SO(n)$ is the
projection. 

The second column of the tables contains the Euler characteristic.

\renewcommand\arraystretch{1.2}
\begin{itemize}
\item[$\bullet$] $\G=\SU(3)$
 \begin{center} \begin{tabular}{|c|c|}\hline
    $\S^1\subset\SO(3),\S(\U(1)\U(2))$ & 3\\\hline
    $\S^1\subset\SO(3),\T^2$ & 6\\\hline
    $\Z_3\SO(2)\subset\Z_3\SO(3),\T^2$ & 6\\\hline
  \end{tabular}\end{center}
  
 \item[$\bullet$] $\G=\SU(4)$
  \begin{center} \begin{tabular}{|c|c|}\hline
    $\S^1\SU(2)\subset\S(\U(2)\U(2)),\S(\U(1)\U(3))$ & 10\\
    where $\S^1=\{\diag(\bar z^2,z^4,\bar z,\bar z)\}\subset N(\SU(2))$ &\\\hline
    $\S^1\subset\sigma(\S(\U(1)\U(3))),\S(\U(1)\U(3))$ & 8\\
    $\S^1=\{\diag(\bar z,z,1,1)\}$, $\sigma$ exchanges the first two coordinates
    & \\\hline
    \end{tabular}\end{center}      
\newpage
\item[$\bullet$] $\G=\SU(n),n\geq 5$
  \begin{center} \begin{tabular}{|c|c|}\hline
    $\S^1\SU(n-2)\subset\sigma(\U(n-1)),\U(n-1)$ & $2n$\\ 
    $\S^1=\{\diag(\bar z,z,1,\dots,1)\}$, $\sigma$ exchanges first two coordinates &\\\hline
    $\S^1\SU(n-2)\subset\S(\U(2)\U(n-2)),\U(n-1)$ & $\frac{n(n+1)}{2}$\\
    $\S^1=\{\diag(\bar z^{n-2},z^{n-2},\bar z,\dots,\bar z)\}\subset N(\SU(n-2))$ & \\\hline
  \end{tabular}\end{center} 
\item[$\bullet$] $\G=\SO(2n+1),n\geq 3$
 \begin{center} \begin{tabular}{|c|c|}\hline
    $\SU(3)\subset\G_2,\U(3),n=3$ & 8 \\\hline
    $\S^1\SU(3)\subset\S^1\G_2,\U(4),n=4$ & 16\\\hline
    $\SO(2)\SO(2n-3)\subset\U(2)\SO(2n-3),\SO(2)\SO(2n-2)$ & $2n(n+1)$\\\hline
    $\T^2\SU(n-2)\subset\SO(3)\S^1\SU(n-2),\SO(2)\U(n-1)$ & $n2^n$\\
    $\S^1=\{\diag(1,1,1,z,\dots,z)\}, \T^2=\{\diag(1,z,z^2,1,\dots,1)\}\cdot
    \S^1$ & \\
    $\SO(3)\into\SO(5)$ is irreducible ($\SO(5)$ upper left block) & \\\hline
 \end{tabular}\end{center} 

\item[$\bullet$] $\G=\Spin(2n+1),n\geq 3$
 \begin{center} \begin{tabular}{|c|c|}\hline
	$\S^1_k\SU(2)\subset\U(3),\SO(2)\SO(5), n=3$ & 14\\
	where $\S^1_k=\{\diag(z^2,z^k,z^k,1)\}$, $k=1,-3$ and &\\
	$\SU(2)$ acts trivially on $\R e_1,\R e_2,\R e_7$ &\\\hline
	$\SU(3)\subset\G_2,\SO(6), n=3$ &2\\\hline 
	$\SU(3)\subset\G_2,\U(3), n=3$ & 8\\\hline
	$\SU(2)\S^1_\H\SU(2)\subset\U(2)\SO(5),\SO(5)\U(2),n=4$& 48\\
	where $\S^1_\H=\{\diag(z^{l_1},z^{l_1},1,z^{l_2},z^{l_2})\}$ and
	$_1=2,l_2=1$ or $l_1=1=l_2$.&\\\hline
 \end{tabular}\end{center} 
 \newpage
\item[$\bullet$] $\G=\Sp(n),n\geq 2$
 \begin{center} \begin{tabular}{|c|c|}\hline
    $\Sp(n-2)\Delta\Sp(1)\subset\Sp(n-2)\SO(2)\Delta\Sp(1),\Sp(n-1)\Sp(1)$ &
    $n(2n-1)$\\
    $\Delta\Sp(1)$ is the diagonal $\Sp(1)$ in the upper left
    $\Sp(2)$-block&\\
    $\SO(2)$ is the standard $\SO(2)$ in this block&\\\hline
    $\Z_2\Delta\Sp(1)\subset\SO(2)\Delta\Sp(1),\Z_2\Sp(1)\Sp(1), n=2$ &8\\
    $\Z_2$ is generated by $\begin{pmatrix}0 & 1\\ 1 &0\end{pmatrix}$
      and &\\
    $\SO(2)$ and $\Delta\Sp(1)$ are as above &\\\hline
    $\S^1\Sp(n-2)\subset\U(2)\Sp(n-2),\S^1\Sp(n-1)$& $2n^2$\\
    where $\S^1=\{\diag(z,\bar z^2,1,\dots,1)\}$ &\\\hline
    $\S^1_\H\Sp(n-2)\subset\Sp(1)\Sp(n-2),\S^1\Sp(n-1)$ & $3n$\\
    $\S^1_\H=\{\diag(z,z^3,1,\dots,1)\}$, $\Sp(1)\into\Sp(2)$ (upper left
    block) &\\
    given by the irreducible $\SO(3)\into \SO(5)$&\\\hline
    $\T^2_\H\Sp(n-3)\subset\S^1\SO'(3)\Sp(n-3),\T^2\Sp(n-2)$ & $4n(n-1)$\\     
    $\S^1=\{\diag(z,z,z,1,\dots,1)\}$, $\T^2_\H=\{\diag(1,z,\bar
    z,1,\dots,1)\}\cdot \S^1$& \\
    $\SO'(3)$ is the conjugation of the standard upper left &\\$\SO(3)$-block
    by $\begin{pmatrix}1 & 0 & 0 \\ 0 & 1 & i\\0& 1&-i\end{pmatrix}$&\\\hline
      \end{tabular}\end{center}    
\item[$\bullet$] $\G=\SO(2n),n\geq 4$
 \begin{center} \begin{tabular}{|c|c|}\hline
   $\SU(4),\Spin(7),\U(4),n=4$ & 8 \\\hline
   $\S^1\SU(4),\S^1\Spin(7),\U(5),n=5$ & 16\\\hline
   $\S^1\SU(n-1)\subset\sigma(\U(n)),\U(n)$ & $2^{n+1}$\\
    where $\S^1=\{\diag(z,1,\dots,1)\}$ or $\S^1=\{\diag(1,z,\dots,z)\}$ and &\\
    $\sigma$ is conjugation by $\diag(-1,1,\dots,1)$&\\\hline
    $\T^2\SU(n-2)\subset\SO(3)\S^1\SU(n-2),\SO(2)\U(n-1)$&$n2^{n-1}$\\
    where $\S^1=\{\diag(1,1,1,1,z,\dots,z)\}$, 
    $\T^2=\{\diag(z,1,\dots,1)\}\cdot \S^1$&\\
    and $\SO(3)$ is the upper left block&\\\hline
  \end{tabular} \end{center}
\end{itemize}
\subsection{Tables} \renewcommand\arraystretch{1.2}
 The following tables contain two well-known classification results: The
 subgroups of maximal rank of the classical lie groups due to Borel and
 Siebenthal, and the transitive effective actions on spheres.
\begin{table}[ht] \begin{center} \begin{tabular}{|l|l|}\hline 
$\G$ &  $\K_+$ \\\hline 
$\SU(n)$ & $\S(\U(n_1)\U(n_2)\cdots \U(n_k))$ \\
&\quad where $n_1+n_2+\ldots+n_k=n$\\\hline 
$\SO(2n)$ & $\SO(2n_1)\cdots\SO(2n_k)\U(m_1)\cdots\U(m_l)$\\ 
&\quad where $\sum n_i + \sum m_i=n$ \\\hline 
$\SO(2n+1)$& $\SO(2k+1)\SO(2n_1)\cdots\SO(2n_k)\U(m_1)\cdots\U(m_l)$\\
&\quad where $k + \sum n_i + \sum m_i=n$\\\hline
$\Sp(n)$ & $\U(n_1)\cdots \U(n_k)\Sp(m_1)\cdots\Sp(m_l)$\\
&\quad where $\sum n_i + \sum m_i=m$\\\hline 
\end{tabular}
\caption{Connected subgroups $\K_+$ of maximal rank of the simple classical Lie
groups}
\label{table:maxranksubgroups}                             
 \end{center}
\end{table}

\begin{table}[ht]
\begin{center}
\begin{tabular}{|c|c|c|c|}\hline
\textbf{Dimension} & \textbf{K} & \textbf{H} & \textbf{Isotropy
representation}\\\hline
n & $\SO(n+1)$ & $\SO(n)$ & $\rho_n$\\\hline
2n+1 & $\SU(n+1)$ & $\SU(n)$ & $\mu_n\oplus id$\\\hline 
2n+1 & $\U(n+1)$ & $\U(n)$ & $\mu_n\oplus id$\\\hline
4n+3 & $\Sp(n+1)$ & $\Sp(n)$ & $\nu_n\oplus 3id$\\\hline
4n+3 & $\Sp(n+1)\Sp(1)$ & $\Sp(n)\Delta\Sp(1)$ & $\nu_n\hat\otimes\nu_1\oplus
id\hat\oplus\rho_3$\\\hline
4n+3 & $\Sp(n+1)\U(1)$& $\Sp(n)\Delta\U(1)$ & $\nu_n\hat\otimes\phi \oplus
id\hat\otimes\phi\oplus id$\\\hline
15 & $\Spin(9)$ & $\Spin(7)$ & $\rho_7\oplus\Delta_8$\\\hline
7 & $\Spin(7)$ & $\G_2$ & $\phi_7$\\\hline
6 & $\G_2$ & $\SU(3)$ & $\mu_3$\\\hline
\end{tabular}\end{center}
\caption{Transitive effective actions on the sphere $S^n$}
\label{table:sphereactions}
\end{table}

\renewcommand\arraystretch{1.0}

\subsection{Actions of $\U(n)$ on spheres}\label{subsubsec:un1spheres}
Actions of $\U(n)$ on spheres of dimension greater  than 1 are given in the
following way: $A\in\U(n)$ acts by $(\det A)^kA$. It is easily seen that the
isotropy group of such an action is given by

\begin{displaymath}
\left\{
  \begin{pmatrix}
    a &0&0&\cdots\\
    0 & & & \\
    0 & & A &\\
    \vdots & & &
  \end{pmatrix}
\mid A\in\C, A\in\U(n-1),a^{k+1}=\overline{\det A}^k \right\}
\end{displaymath}

Also one easily verifies that this is $\S^1_k\SU(n-1)$ where

\begin{displaymath}
  \S^1_k=\{\diag(z^{(k+1)(n-1)},z^{-k},\dots,z^{-k})\mid z\in\S^1\}
\end{displaymath}

Replace $z$ by $z^n$ and write this as a product of

\begin{displaymath}
  \diag(z^{n-1},z^{n-1},\dots,z^{n-1})
\end{displaymath}

and

\begin{displaymath}
  \diag(z^{(k+1)(n-1)n-(n-1)},
  z^{-kn-(n-1)},\dots,z^{-kn-(n-1)})
\end{displaymath}

where the first is an element of the center and the second in $\SU(n)$.

\subsection{Actions of $\U(n)$ on real projective
spaces}\label{subsubsec:u(n)onproj}

If $\U(n)$ acts transitively almost effectively on the real projective space $\R
P^{2n-1}$, the isotropy group is given by $\S^1_{l,k}\SU(n-1)$, where

\begin{displaymath}
  \S^1_{l,k}=\{\diag(z^{l},z^{-k},\dots,z^{-k})\mid z\in\S^1\},\gcd(k,l)=1
\end{displaymath}

We want to show $l=(k+2)(n-1)$. 

To prove that, we look at the map $\pi_1(\S^1_{l,k})\to\pi_1(\U(n))$ induced by
inclusion. We have $\pi_1(\R P^{2n-1})=\Z_2$, so the image has index 2 in the
fundamental group of $\U(n)$. We know $\S^1_{l,k}$ is homotopic to

\begin{displaymath}
  \{\diag(z^{l},z^{-k(n-1)},1\dots,1)\}=\{\diag(z^{\frac{l}{a}},z^{-k\frac{n-1}{a}},1,\dots,1)\}
\end{displaymath}

where $a=\gcd(l,n-1)$. This in turn is homotopic to

\begin{displaymath}
  \{\diag(z^{\frac{l}{a}-k\frac{n-1}{a}},1,\dots,1)\}
\end{displaymath}

By what we said above, it is clear that $|\frac{l}{a}-k\frac{n-1}{a}|=2$. First
suppose $\frac{l}{a}-k\frac{n-1}{a}=2$. Then $\frac la=k\frac{n-1}a+2$, and we
consider the following sequence (where homotopy is denoted by $\simeq$):

\begin{align*}
  &\{\diag(z^{k\frac{n-1}a},z^{-k\frac{n-1}a},1,\dots,1)\}\\
  =&\{\diag(z^{(k\frac{n-1}a+2)(n-1)},z^{-k(n-1)\frac{n-1}a},1,\dots,1)\}\\
  \simeq
  &\{\diag(z^{(k\frac{n-1}a+2)(n-1)},z^{-k\frac{n-1}a},z^{-k\frac{n-1}a},\dots,z^{-k\frac{n-1}a})\}
\end{align*}

By replacing $k$ by $k\frac{n-1}a$, we arrive at the aforementioned form. A
similar computation shows the same if $\frac{l}{a}-k\frac{n-1}{a}=-2$.

The converse of the claim is also clear from the above computation: The index of
the image of $\S^1_{l,k}$ in the fundamental group of $\U(n)$ has index 2, if
$l=(k+2)(n-1)$ and $\gcd(l,k)=1$.

\bibliographystyle{plainnat}
\bibliography{cohom1}
\end{document}